\RequirePackage{fix-cm}
\documentclass[11pt,twoside,english]{article}
\usepackage[T1]{fontenc}
\usepackage[utf8]{inputenc}
\usepackage[a4paper]{geometry}
\newcommand{\coshh}{ \cosh \left( \sqrt{\lambda_p}t\right)}
\newcommand{\sinhh}{ \sinh \left( \sqrt{\lambda_p}t\right)}

\newcommand{\sinhs}{ \sinh \left( \sqrt{\lambda_p}(t-s)\right)}
\newcommand{\la}{\lambda}
\newcommand{\canla}{  \sqrt{\lambda_p} }
\newcommand{\be}{\begin{eqnarray*}}                           
\newcommand{\en}{\end{eqnarray*}}
\newcommand{\bes}{\begin{eqnarray}}                           
\newcommand{\ens}{\end{eqnarray}}
\usepackage{babel}
\usepackage{array}
\usepackage{endnotes}
\usepackage{multirow}
\usepackage{amsthm}
\usepackage{amsmath}
\usepackage{amssymb}
\usepackage{fixltx2e}
\usepackage{graphicx}
\usepackage{esint}
\usepackage[unicode=true,pdfusetitle,
 bookmarks=true,bookmarksnumbered=true,bookmarksopen=true,bookmarksopenlevel=1,
 breaklinks=false,pdfborder={0 0 1},backref=false,colorlinks=false]
 {hyperref}
\usepackage{breakurl}
\makeatletter

\providecommand{\tabularnewline}{\\}

\usepackage{enumitem}		
\usepackage{braille}

\@ifundefined{lettrine}{\usepackage{lettrine}}{}
 \let\footnote=\endnote
\theoremstyle{plain}
\ifx\thechapter\undefined
\newtheorem{thm}{\protect\theoremname}
\else
\newtheorem{thm}{\protect\theoremname}[chapter]
\fi
  \theoremstyle{remark}
  \newtheorem{rem}[thm]{\protect\remarkname}
  \theoremstyle{plain}
  \newtheorem{lem}[thm]{\protect\lemmaname}

\makeatother

  \providecommand{\lemmaname}{Lemma}
  \providecommand{\remarkname}{Remark}
\providecommand{\theoremname}{Theorem}

\begin{document}

\title{\textbf{Approximation of mild solutions of the linear and nonlinear
elliptic equations}}

\author{Nguyen Huy Tuan\textsuperscript{1,}, Dang Duc Trong\textsuperscript{1},
Le Duc Thang\textsuperscript{2} and Vo Anh Khoa\textsuperscript{1}\\\\
{\small{\textsuperscript{1}Department of Mathematics and Computer
Science, University of Science,}}\\
{\small{227 Nguyen Van Cu Street, District 5, Ho Chi Minh City, Vietnam.}}\\
{\small{\textsuperscript{2}Faculty of Basic Science, Ho Chi Minh
City Industry and Trade College,}}\\
{\small{20 Tang Nhon Phu, District 9, Ho Chi Minh City, Viet Nam.}}}
\maketitle
\begin{abstract}
In this paper, we investigate the Cauchy problem for both linear and
semi-linear elliptic equations. In general, the equations have the
form

\[
\frac{\partial^{2}}{\partial t^{2}}u\left(t\right)=\mathcal{A}u\left(t\right)+f\left(t,u\left(t\right)\right),\quad t\in\left[0,T\right],
\]

where $\mathcal{A}$ is a positive-definite, self-adjoint operator
with compact inverse. As we know, these problems are well-known to
be ill-posed. On account of the orthonormal eigenbasis and the corresponding
eigenvalues related to the operator, the method of separation of variables
is used to show the solution in series representation. Thereby, we
propose a modified method and show error estimations in many accepted
cases. For illustration, two numerical examples, a modified Helmholtz
equation and an elliptic sine-Gordon equation, are constructed to
demonstrate the feasibility and efficiency of the proposed method.


\end{abstract}

\section{Introduction}

The Cauchy problem of elliptic equation plays an important role in
inverse problems. For example, in optoelectronics, the determination
of a radiation field surrounding a source of radiation (e.g., a light
emitting diode) is a frequently occurring problem. As a rule, experimental
determination of the whole radiation field is not possible. Practically,
we are able to measure the electromagnetic field only on some subset
of physical space (e.g., on some surfaces). So, the problem arises
how to reconstruct the radiation field from such experimental data
(see, for instance, \cite{key-22}). In the paper of Reginska \cite{key-22},
the authors considered a physical problem which is connected with
the notion of light beams. Some applications of this model can be
established in more detail in \cite{key-22}. Another application
in inverse obstacle problems (cf. \cite{key-4}), which are investigated
in connection with inclusion detection by electrical impedance tomography
when only one pair of boundary current and voltage is used for probing
the examined body \cite{key-20}. 

Let $\mathcal{H}$ be a real Hilbert space, and let $\mathcal{A}:\mathcal{D}\left(\mathcal{A}\right)\subset\mathcal{H}\to\mathcal{H}$
be a positive-definite, self-adjoint operator with compact inverse
on $\mathcal{H}$. In this paper, we consider the problem of finding
a function $u:\left[0,T\right]\to\mathcal{H}$ satisfying

\begin{equation}
\frac{\partial^{2}}{\partial t^{2}}u\left(t\right)=\mathcal{A}u\left(t\right)+f\left(t,u\left(t\right)\right),\quad t\in\left[0,T\right],\label{eq:1}
\end{equation}

associated with the initial conditions

\begin{equation}
u\left(0\right)=\varphi,\quad\frac{\partial}{\partial t}u\left(0\right)=g,\label{eq:2}
\end{equation}

where $f$ is a mapping from $\left[0,T\right]\times\mathcal{H}\to\mathcal{H}$,
$\varphi$ and $g$ are the exact data in $\mathcal{H}$. Physically,
the exact data can only be measured, there will be measurement errors,
and we thus would have as data some function $\varphi^{\epsilon}$
and $g^{\epsilon}$ in $\mathcal{H}$ for which

\begin{equation}
\left\Vert \varphi-\varphi^{\epsilon}\right\Vert \le\epsilon,\quad\left\Vert g-g^{\epsilon}\right\Vert \le\epsilon,\label{eq:3}
\end{equation}

where the constant $\epsilon>0$ represents a bound on the measurement
error, $\left\Vert .\right\Vert $ denotes the $\mathcal{H}$ norm.\\

 Since Hadamard\cite{Ha}, it is well known that the Cauchy problem of elliptic equation, for example, Problem (1)-(2), 
 is severely ill-posed: although it has at most one solution, it may
 have none, and if a solution exists, it does not depend continuously on the data
$\varphi, g$ in any reasonable topology. Therefore, regularization is needed to stabilize
 the problem.
 In recent years, many special regularization methods for the homogeneous
 and nonhomogeneous Cauchy problem of elliptic equation have been proposed,
 such as Backus-Gilbert algorithm \cite{key-9}, the method of wavelet
 \cite{key-12}, quasi-reversibility method \cite{key-17}, truncation
 method \cite{key-25}, non-local boundary value method \cite{key-10}
 and the references therein.

Although we have many works on the linear homogeneous case of Cauchy  problem for elliptic equation, however, regularization theory and numerical
simulation for nonlinear elliptic equations are still limited. Especially,
the nonlinear cases for elliptic equation appear in many real applications.   For example, let us see a simple one infered by giving $\mathcal{A}=\dfrac{-\partial^{2}}{\partial x^{2}}$
and $\mathcal{D}\left(\mathcal{A}\right)=H_{0}^{1}\left(0,\pi\right)\subset\mathcal{H}=L^{2}\left(0,\pi\right)$
in the problem (\ref{eq:1})-(\ref{eq:2}). In particular, it is given
by

\begin{equation}
\begin{cases}
\frac{\partial^{2}}{\partial t^{2}}u\left(x,t\right)+\frac{\partial^{2}}{\partial x^{2}}u\left(x,t\right)=f\left(x,t,u\left(x,t\right)\right) & ,\left(x,t\right)\in\left(0,\pi\right)\times\left(0,1\right),\\
u\left(0,t\right)=u\left(\pi,t\right)=0 & ,t\in\left(0,1\right),\\
u\left(x,0\right)=\varphi\left(x\right),\quad\frac{\partial}{\partial t}u\left(x,0\right)=g\left(x\right) & ,x\in\left(0,\pi\right).
\end{cases}\label{eq:4}
\end{equation}
If $f\left(x,t,u\right)=k^{2}u$ in (\ref{eq:4}), then it
is called Helmholtz equation which has many applications related to
wave propagation and vibration phenomena. This equation is often used
to describe the vibration of a structure, the acoustic cavity problem,
the radiated wave and the scattering of a wave. With $f\left(x,t,u\right)=\sin u$ in (\ref{eq:4}), we obtain the
elliptic sine-Gordon equation. From the point of view of the modelling
of physical phenomena, the motivation for the study of this equation
comes from its applications in several areas of mathematical physics
including the theory of Josephson effects, superconductors and spin
waves in ferromagnets, see e.g. \cite{key-15}. With $f\left(x,t,u\right)=u-u^{3}$,
we have the Allen-Cahn equation originally formulated in the description
of bi-phase separation in fluids.

Switch back to the considered problem, it is more complicated than
the ones above. Hence, the purpose of this paper is to introduce a
new method of integral equation that is based on a modification of
the exact solution formulation. As the regularization parameter tends
to zero, the solution of our regularized problem converges monotonically
to the solution of the Cauchy problem with the exact data.

Prior to the approach of main results, we would like to introduce
the representation of solution in problem (\ref{eq:1})-(\ref{eq:2})
for linear and semi-linear cases. We can see that the operator $\mathcal{A}$,
as a consequence, admits an orthonormal eigenbasis $\left\{ \phi_{p}\right\} _{p\ge1}$
in $\mathcal{H}$, associated with the eigenvalues such that

\begin{equation}
0<\lambda_{1}\le\lambda_{2}\le...\lim_{p\to\infty}\lambda_{p}=\infty.
\end{equation}
Let $u(t)= \sum\limits_{p=1}^\infty \left\langle u(t),\phi_p\right\rangle \phi_p$ be the Fourier series of $u$ in the Hilbert space $H$. For homogeneous problem, i.e, $f=0$ in (1), by a seperable method,  we get the homogeneous second order differential equation as follows
\be
\frac{d^2}{dt^2} \left\langle u(t),\phi_p \right\rangle -\la_p \left\langle u(t),\phi_p \right\rangle = 0, \left\langle u(0),\phi_p \right\rangle= \left\langle \varphi,\phi_p \right\rangle,~\frac{d}{dt} \left\langle u(0),\phi_p \right\rangle=  \left\langle g,\phi_p \right\rangle,
\en    
and its solution leads to 
\begin{eqnarray} 
u(t)=\sum\limits_{p=1}^\infty \left[\coshh \left\langle \varphi,\phi_p \right\rangle + \frac{\sinhh}{\canla} \left\langle g,\phi_p \right\rangle\right]\phi_p, \label{eq:6}
\end{eqnarray} 
where $\left\langle .,.\right\rangle $ denotes the inner product
in $\mathcal{H}$.
From F. Browder terminology, as in {[}Dan Henry,\emph{ Geometric Theory
of Semi-linear Parabolic Equations}, Springer-Verlag, Berlin Heildellberg,
Berlin, 1982{]}, $u(t)$  in \eqref{eq:6} is called the \emph{mild solution}
of (\ref{eq:1})-(\ref{eq:2}) with $f=0$.

For the nonlinear problem $(1)-(2)$, we say that $ u \in C([0,T];H) $ is a mild solution  if $u$ satisfies the integral equation 
\begin{eqnarray} 
&& u(t)=\sum\limits_{p=1}^\infty \left[\coshh \varphi_p + \frac{\sinhh}{\canla} g_p+\int\limits_{0}^{t}\frac{\sinhs}{\canla} f_{p}(u)(s)ds\right]\phi_p \label{eq:7}
\end{eqnarray}
where $f_{p}(u)(s) = \left\langle f(s,u(s)), \phi_p)\right\rangle$. 
The transformation  from problem (1)-(2) into \eqref{eq:7} is easily proved by a separation method which is similar above process. From now on, to regularize Problem (1)-(2), we only consider the integral equation \eqref{eq:7} and find a regularization method for it. The main idea of integral equation method can be found in a paper \cite{Trong} on nonlinear backward heat equation.

The paper is organized as follows. In Section 2, we present our regularization
method for the linear problem implied by letting $f=0$ in (\ref{eq:1}).
The theoretical results in the Section 2 are inspirable for us to
suggest a new regularization method for semi-linear case in Section
3. New convergence estimates are given under some different priori
assumptions for the exact solution. Proofs of the results in these
sections will be showed in the appendix in the bottom of paper. In
Section 4, simple numerical examples aimed to illustrate the main
results in Section 3 are analyzed.

\section{The linear homogeneous problem}
In \cite{key-17}, C.L. Fu and his group applied the {\it quasi-reversibiity } (QR ) method to approximate problem (4) in case $f=0$ and $g=0$. The main idea of the original QR method \cite{lions} is to approach the ill-posed
second order Cauchy problem by a family of well-posed fourth order problems depending
on a (small) regularization parameter. In particular, they considered approximate problem

\begin{equation}
\begin{cases}
u^\epsilon_{tt}\left(x,t\right)+ u^\epsilon_{xx}\left(x,t\right)-\beta^2 u^\epsilon_{ttxx}\left(x,t\right)  =0 & ,\left(x,t\right)\in\left(0,\pi\right)\times\left(0,1\right),\\
u\left(0,t\right)=u\left(\pi,t\right)=0 & ,t\in\left(0,1\right),\\
u\left(x,0\right)=\varphi^\epsilon\left(x\right),\quad\frac{\partial}{\partial t}u\left(x,0\right)=0 & ,x\in\left(0,\pi\right).
\end{cases}\label{eq:444}
\end{equation}
The solution of \eqref{eq:444} is defined by
\bes
u^\epsilon(x,t)=  \sum_{p=1}^\infty \cosh \left( \frac{pt}{\sqrt{1+\beta^2 p^2}} \right)  \left<\varphi^\epsilon(x), \sin (px)\right> \sin(px)  \label{fu}
\ens
and the authors proved that $u^\epsilon$ converges to the solution $u$ of homogeneous problem as $\epsilon \to 0$.\\
Very recently, homogeneous problem has been considered by Hao, Duc and  Lesnic  \cite{key-10}. They applied the method of {\it non-local boundary value problems} (also called quasi-boundary value method) to regularized the above problem as follows
\bes
\left\{ \begin{gathered}
  u_{tt}  =Au, \hfill \\
  u_t(0)=0 \hfill \\
 u(0)+ \beta u(aT)=\varphi\hfill\\
 \end{gathered}  \right.\label{e9}
\ens
with $a \ge 1 $ being given and $\beta >0$ is the regularization parameter. They proved that the solution to \eqref{e9} is
\bes
u^\alpha(t)=  \sum_{p=1}^\infty  \frac{\coshh}{1+\beta \cosh (a \sqrt{\lambda_p}t)}  \Big<\varphi, \phi_p \Big> \phi_p
\ens
and $\|u^\beta(t)-u(t)\| \to 0$ as $\beta  \to 0$ with some assumptions on the exact solution $u$. \\
Following the work \cite{key-10}, in \cite{key-25} Tuan, Trong and Quan used a Fourier truncated method to
treat the following Cauchy problem of an elliptic equation with nonhomogeneous
Dirichlet and Neumann data. From the simple analysis about the exact solution \eqref{eq:6}, we know that the data error can be arbitrarily amplified by the “kernel” $\coshh$. That is the reason
why the Cauchy problem of elliptic equation is ill-posed. Since the general regularization theory \cite{ki} and paper \cite{key-17}, we now give a more general principle
of regularization methods for the Cauchy problem of (6). Our idea on
regularization method is of  constructing a new kernel $Q(t,\lambda_p,\beta)$ and replacing $\coshh$ by  $Q(t,\lambda_p,\beta)$ where the new
kernel should satisfy\\

(A) If $\beta$ is fixed, $  Q(t,\lambda_p,\beta) $ is bounded.\\
(B)  If $t, \lambda_p$ is fixed, then $\lim_{\beta \to 0}  Q(t,\lambda_p,\beta)=\coshh  $.\\

Following properties (A) and (B), one can construct other kernels. Furthermore, the idea of properties (A) and
(B) can be applied to other ill-posed problems when the solution has the similar form of (6), e.g., the inverse heat
conduction problem \cite{Qian2}. In this sense, we say that the properties (A) and (B) are useful and interesting. Now, from above discussion, it is easy to check   that the kernels $Q_1(t,\lambda_p,\beta)= \cosh \left( \frac{\sqrt{\lambda_p} t}{\sqrt{1+\beta^2 \lambda_p}} \right)$ in \cite{key-17} and $Q_2(t,\lambda_p,\beta)= \frac{\coshh}{1+\beta \cosh (a \sqrt{\lambda_p}t)} $ in \cite{key-10} satisfy (A) and (B).\\

We now have a look at the solution $u$ in (\ref{eq:6}). To find a  regularization solution for $u$, the unstability terms $\coshh$ and $\sinhh$ in (6) should be  replaced by two kernels $Q(t,\lambda_p,\beta)$ and $ R(t,\lambda_p,\beta)$  respectively. Here the kernel $Q$ satisfies (A), (B) and kernel $R$ satisfies the following conditions

(C) If $\beta$ is fixed, $  R(t,\lambda_p,\beta) $ is bounded.\\
(D)  If $t, \lambda_p$ is fixed, then $\lim_{\beta \to 0}  R(t,\lambda_p,\beta)=\sinhh  $.\\

In \cite{key-25}, we choose 
\bes
Q(t,\lambda_p,\beta) = R(t,\lambda_p,\beta)= 
\left\{ \begin{gathered}
 1,~~\text {if}~~\la_p \le m_\beta^2, \hfill \\
0,~~\text {if}~~\la_p > m_\beta^2, \hfill\\
 \end{gathered}  \right. \label{filter2}
\ens
to get a truncation solution (See the fomula  (7) in  page 2915, \cite{key-25} ) where $m_\beta$ such that $\lim_{\beta \to +\infty m_\beta }=+\infty$. It is easy to check that $Q$ and $R$ defined in \eqref{filter2} satisfy $(A), (B)$ and $(C), (D)$ respectively.  \\

In this section, we consider the homogeneous problem of (1) (also  given in \cite{key-25} ) by other choices for kernels. From the formula of $\coshh$ and $\sinhh$, we realize
that the term $e^{\sqrt{\lambda_{p}}t}$ is unstability cause while
the term $e^{-\sqrt{\lambda_{p}}t}$ is stable under the boundedness
of the unity. Hence, by a simple and natural way, we  replace $\coshh$ and $\sinhh$ by two new kernels 
\be
Q_3(t,\lambda_p,\beta)= \dfrac{1}{2\beta+2e^{-\sqrt{\lambda_{p}}t}}+ \frac{e^{-\sqrt{\lambda_{p}}t}}{2},
\en
and 
\be
R(t,\lambda_p,\beta)= \dfrac{1}{2\beta+2e^{-\sqrt{\lambda_{p}}t}}-\frac{e^{-\sqrt{\lambda_{p}}t}}{2},
\en
to obtain a regularization solution
\begin{eqnarray} 
u^\epsilon(t)=\sum\limits_{p=1}^\infty \left[ Q_3(t,\lambda_p,\beta) \left\langle \varphi,\phi_p \right\rangle + \frac{R(t,\lambda_p,\beta)}{\canla} \left\langle g,\phi_p \right\rangle\right]\phi_p. \label{resol}
\end{eqnarray} 
Here $\beta=\beta\left(\epsilon\right)$ is called parameter reguarization and satisfies ${\displaystyle \lim_{\epsilon\to0}\beta\left(\epsilon\right)=0}$. It is easy to check that $Q_3$ and $R$ satisfy $(A), (B)$ and $(C), (D)$ respectively.  Moreover, \eqref{resol} leads to

\begin{equation}
u^{\epsilon}\left(t\right)=\sum_{p\ge1}\left[\dfrac{1}{2\beta+2e^{-\sqrt{\lambda_{p}}t}}\left(\left\langle \varphi,\phi_{p}\right\rangle +\frac{\left\langle g,\phi_{p}\right\rangle }{\sqrt{\lambda_{p}}}\right)+\frac{e^{-\sqrt{\lambda_{p}}t}}{2}\left(\left\langle \varphi,\phi_{p}\right\rangle -\frac{\left\langle g,\phi_{p}\right\rangle }{\sqrt{\lambda_{p}}}\right)\right]\phi_{p}.\label{eq:8}
\end{equation}

Under the inexact data $\varphi^{\epsilon}$ and $g^{\epsilon}$,
the regularized solution becomes

\begin{equation}
v^{\epsilon}\left(t\right)=\sum_{p\ge1}\left[\dfrac{1}{2\beta+2e^{-\sqrt{\lambda_{p}}t}}\left(\left\langle \varphi^{\epsilon},\phi_{p}\right\rangle +\frac{\left\langle g^{\epsilon},\phi_{p}\right\rangle }{\sqrt{\lambda_{p}}}\right)+\frac{e^{-\sqrt{\lambda_{p}}t}}{2}\left(\left\langle \varphi^{\epsilon},\phi_{p}\right\rangle -\frac{\left\langle g^{\epsilon},\phi_{p}\right\rangle }{\sqrt{\lambda_{p}}}\right)\right]\phi_{p}.\label{eq:9}
\end{equation}

\begin{rem}
\label{rm1}With this linear case of (\ref{eq:1}) we denote the solution
of (\ref{eq:1})-(\ref{eq:2}) by $u\left(t\right)$, the regularized
solution of (\ref{eq:1})-(\ref{eq:2}) by $u^{\epsilon}\left(t\right)$,
and the regularized solution of (\ref{eq:1})-(\ref{eq:3}) by $v^{\epsilon}\left(t\right)$.
\end{rem}
The main results of this section are in the following theorem.
\begin{thm}
\label{thm:2}Let $\beta=\epsilon^{m}$ for $m\in\left(0,1\right)$.\end{thm}
\begin{description}
\item [{(i)}] If there is a positive constant $E_{1}$ such that

\begin{equation}
\sqrt{\frac{\left\Vert u\left(T\right)\right\Vert ^{2}}{2}+\frac{\left\Vert \dfrac{\partial}{\partial t}u\left(T\right)\right\Vert ^{2}}{2\lambda_{1}}}<E_{1},\label{eq:10}
\end{equation}

then we have

\begin{equation}
\begin{cases}
\left\Vert u\left(t\right)-v^{\epsilon}\left(t\right)\right\Vert \le\sqrt{2\left(1+\frac{1}{\lambda_{1}}\right)}\epsilon^{1-m}+E_{1}\epsilon^{m} & ,t\in\left[0,\frac{T}{2}\right],\\
\left\Vert u\left(t\right)-v^{\epsilon}\left(t\right)\right\Vert \le\sqrt{2\left(1+\frac{1}{\lambda_{1}}\right)}\epsilon^{1-m}+E_{1}\epsilon^{\frac{m\left(T-t\right)}{t}} & ,t\in\left[\frac{T}{2},T\right].
\end{cases}
\end{equation}

\item [{(ii)}] If there is a positive constant $E_{2}$ such that

\begin{equation}
\sqrt{\sum_{p\ge1}e^{2\sqrt{\lambda_{p}}\left(T-t\right)}\left(\sqrt{\lambda_{p}}\left\langle u\left(t\right),\phi_{p}\right\rangle +\left\langle \dfrac{\partial}{\partial t}u\left(t\right),\phi_{p}\right\rangle \right)^{2}}<E_{2},\label{eq:12}
\end{equation}

then we have

\begin{equation}
\begin{cases}
\left\Vert u\left(t\right)-v^{\epsilon}\left(t\right)\right\Vert \le\sqrt{2\left(1+\frac{1}{\lambda_{1}}\right)}\epsilon^{1-m}+\frac{\epsilon^{m}}{2\sqrt{\lambda_{1}}}E_{2} & ,t\in\left[0,\frac{T}{2}\right],\\
\left\Vert u\left(t\right)-v^{\epsilon}\left(t\right)\right\Vert \le\sqrt{2\left(1+\frac{1}{\lambda_{1}}\right)}\epsilon^{1-m}+\frac{\epsilon^{\frac{m\left(T-t\right)}{t}}}{2\sqrt{\lambda_{1}}}\left[\frac{\lambda_{1}T}{1+\ln\left(\frac{\sqrt{\lambda_{1}}T}{\epsilon^{m}}\right)}\right]^{\frac{2t-T}{t}}E_{2} & ,t\in\left[\frac{T}{2},T\right].
\end{cases}
\end{equation}

\item [{(iii)}] If there is a positive constant $E_{3}$ such that

\begin{equation}
\sqrt{\sum_{p\ge1}e^{2\sqrt{\lambda_{p}}t}\left(\left\langle \sqrt{\lambda_{p}}u\left(t\right)+\dfrac{\partial}{\partial t}u\left(t\right),\phi_{p}\right\rangle \right)^{2}}<E_{3},\label{eq:14}
\end{equation}

then we have

\begin{equation}
\left\Vert u\left(t\right)-v^{\epsilon}\left(t\right)\right\Vert \le\sqrt{2\left(1+\frac{1}{\lambda_{1}}\right)}\epsilon^{1-m}+E_{3}\frac{\epsilon^{m}}{2}.\label{eq:15}
\end{equation}

\end{description}
In order to prove this theorem, we have to obtain some auxiliary results
given by the lemmas below.
\begin{lem}
\label{lem:3}Let $0<\beta<1$ and let $u^{\epsilon}\left(t\right),v^{\epsilon}\left(t\right)\in\mathcal{H}$
as introduced in Remark \ref{rm1}. Then, we have the following estimate

\begin{equation}
\left\Vert u^{\epsilon}\left(t\right)-v^{\epsilon}\left(t\right)\right\Vert \le\sqrt{2\left(1+\frac{1}{\lambda_{1}}\right)}\epsilon\beta^{-1}.
\end{equation}

\end{lem}

\begin{lem}
\label{lem:4}Let $0<\beta<1$ and let $u\left(t\right),v^{\epsilon}\left(t\right)\in\mathcal{H}$
as introduced in Remark \ref{rm1}. If (\ref{eq:10}) is satisfied,
then we have the following estimate

\begin{equation}
\begin{cases}
\left\Vert u\left(t\right)-v^{\epsilon}\left(t\right)\right\Vert \le\beta E_{1} & ,t\in\left[0,\frac{T}{2}\right],\\
\left\Vert u\left(t\right)-v^{\epsilon}\left(t\right)\right\Vert \le\beta^{\frac{T-t}{t}}E_{1} & ,t\in\left[\frac{T}{2},T\right].
\end{cases}\label{eq:17}
\end{equation}

\end{lem}

\begin{lem}
\label{lem:5}Let $0<\beta<1$ and let $u\left(t\right),v^{\epsilon}\left(t\right)\in\mathcal{H}$
as introduced in Remark \ref{rm1}. If (\ref{eq:12}) is satisfied,
then we have

\begin{equation}
\begin{cases}
\left\Vert u\left(t\right)-v^{\epsilon}\left(t\right)\right\Vert \le\frac{\beta}{2\sqrt{\lambda_{1}}}E_{2} & ,t\in\left[0,\frac{T}{2}\right],\\
\left\Vert u\left(t\right)-v^{\epsilon}\left(t\right)\right\Vert \le\frac{1}{2\sqrt{\lambda_{1}}}\beta^{\frac{T-t}{t}}E_{2} & ,t\in\left[\frac{T}{2},T\right].
\end{cases}\label{eq:18}
\end{equation}

\end{lem}

\begin{lem}
\label{lem:6}Let $0<\beta<1$ and let $u\left(t\right),v^{\epsilon}\left(t\right)\in\mathcal{H}$
as introduced in Remark \ref{rm1}. If (\ref{eq:14}) is satisfied,
then the following estimate holds

\begin{equation}
\left\Vert u\left(t\right)-v^{\epsilon}\left(t\right)\right\Vert \le\frac{\beta}{2}E_{3}.\label{eq:19}
\end{equation}
\end{lem}
\begin{rem}
At $t=T$, the error in case \textbf{(i) }is useless while it is useful
in case \textbf{(ii)}. Moreover, in case \textbf{(iii)}, under the
strong assumptions of $u$, we get the error of Holder-logarithmic
type. In fact, if $\epsilon$ is fixed then the right-hand side of
(\ref{eq:15}) get its maximum value at $m=\dfrac{1}{2}$. Thus, we
obtain the error of order $\epsilon^{\frac{1}{2}}$.

On the other hand, the condition in (\ref{eq:12}) is accepted and
natural. Thus, we prove that

\begin{equation}
e^{\sqrt{\lambda_{p}}\left(T-t\right)}\left(\sqrt{\lambda_{p}}\left\langle u\left(t\right),\phi_{p}\right\rangle +\left\langle \frac{\partial}{\partial t}u\left(t\right),\phi_{p}\right\rangle \right)=\sqrt{\lambda_{p}}\left\langle u\left(T\right),\phi_{p}\right\rangle +\left\langle \frac{\partial}{\partial t}u\left(T\right),\phi_{p}\right\rangle .
\end{equation}

Then the condition 

\begin{equation}
\sum_{p\ge1}\left(\sqrt{\lambda_{p}}\left\langle u\left(T\right),\phi_{p}\right\rangle +\left\langle \frac{\partial}{\partial t}u\left(T\right),\phi_{p}\right\rangle \right)^{2}<\infty,
\end{equation}

is easy to check.
\end{rem}

\section{The semi-linear problem}

As we introduced, many previous papers only regularized problems related
to (\ref{eq:1}) in which $f=0$. This condition makes the applicability
of the method very narrow. Until now, the results in nonlinear case
are very rare. In this section, we consider the problem (\ref{eq:1})
where $f:R\times\mathcal{H}\to\mathcal{H}$ is a Lipschitz continuous
function, i.e., there exists $K>0$ independent of $w_{1},w_{2}\in\mathcal{H},t\in R$
such that

\begin{equation}
\left\Vert f\left(t,w_{1}\right)-f\left(t,w_{2}\right)\right\Vert \le K\left\Vert w_{1}-w_{2}\right\Vert .
\end{equation}

Since
$0< t < s \le T$ , we know from (\ref{eq:7})   that, when $p$ becomes large, the terms $$\coshh,\sinhh,\sinhs,$$ increases rather quickly. Thus, these terms
are the
unstability causes. Hence, to find a regularization solution, we have to replace these terms by new kernels (called stability terms). These kernels have some common properties $(A), (B), (C), (D)$. In fact, we define a following regularization solution
\begin{eqnarray} 
&& u^\epsilon(t)=\sum_{p\ge1} \left[ P(t,\lambda_p, \beta) \varphi_p + \frac{Q(t,\lambda_p, \beta)}{\canla} g_p+\int\limits_{0}^{t}\frac{R(t,s,\lambda_p, \beta)}{\canla} f_{p}(u^\epsilon)(s)ds\right]\phi_p. \label{resol2}
\end{eqnarray}
Here, $P(t,\lambda_p, \beta),~Q(t,\lambda_p, \beta),~R(t,s,\lambda_p, \beta)$ are bounded by $C(\beta)$ for any $\lambda_p >0$. Moreover,  if $t,~ \lambda_p$ fixed then
\be
&&\lim_{\beta \to 0}  P(t,\lambda_p, \beta) = \coshh,~\lim_{\beta \to 0}  Q(t,\lambda_p, \beta)=\sinhh,\\
&&\lim_{\beta \to 0}R(t,s,\lambda_p, \beta)=\sinhs.
\en
By direct computation, we see that the kernels in Theorem 1 is not applied to nonlinear problem. For solving this problem, we find  some suitable kernels as follows
\be
P(t,\lambda_p, \beta)&=& \frac{e^{-\sqrt{\lambda_{p}}\left(T-t\right)}}{2\beta\sqrt{\lambda_{p}}+2e^{-\sqrt{\lambda_{p}}T}}+\frac{e^{-\sqrt{\lambda_p}t}}{2}, \\
Q(t,\lambda_p, \beta)&=& \frac{e^{-\sqrt{\lambda_{p}}\left(T-t\right)}}{2\beta\sqrt{\lambda_{p}}+2e^{-\sqrt{\lambda_{p}}T}}-\frac{e^{-\sqrt{\lambda_p}t}}{2}, \\
R(t,s,\lambda_p, \beta)&=& \frac{e^{-\sqrt{\lambda_{p}}\left(T+s-t\right)}}{2\beta\sqrt{\lambda_{p}}+2e^{-\sqrt{\lambda_{p}}T}}-\frac{e^{-\sqrt{\lambda_p}(t-s)}}{2}. \\
\en

Then, we show error estimates between the solution $u\left(t\right)$
and the regularized solution $v^{\epsilon}\left(t\right)$ in $\mathcal{H}$
norm under some supplementary error estimates and assumptions. Simultaneously,
the uniqueness of solution $u^{\epsilon},v^{\epsilon}\in C\left(\left[0,T\right];\mathcal{H}\right)$
is proved by contraction principle.

Generally speaking, we obtain the following theorem.
\begin{thm}
\label{thm:8}Let ${\displaystyle u\left(t\right)=\sum_{p\ge1}\left\langle u\left(t\right),\phi_{p}\right\rangle \phi_{p}}$
be the solution as denoted in (\ref{eq:7}). Suppose there is a positive
constant $P$ such that

\begin{equation}
4\sup_{0\le t\le T}\sum_{p\ge1}e^{\sqrt{\lambda_{p}}\left(T-t\right)}\left(\sqrt{\lambda_{p}}\left\langle u\left(t\right),\phi_{p}\right\rangle +\left\langle \frac{\partial}{\partial t}u\left(t\right),\phi_{p}\right\rangle \right)^{2}\le P.\label{eq:23}
\end{equation}

Then by letting $\beta=\epsilon^{m},m\in\left(0,1\right)$ the problem

\begin{eqnarray}
v^{\epsilon}\left(t\right) & = & \sum_{p\ge1}\left[\Phi\left(\beta,\lambda_{p},t\right)\mathcal{M}_{p}\left(\varphi^{\epsilon},g^{\epsilon}\right)+\int_{0}^{t}\Psi\left(\beta,\lambda_{p},s,t\right)\left\langle f\left(s,v^{\epsilon}\left(s\right)\right),\phi_{p}\right\rangle ds\right]\phi_{p}\nonumber \\
 &  & +\sum_{p\ge1}\left[\frac{e^{-\sqrt{\lambda_{p}}t}}{2}\mathcal{M}_{p}\left(\varphi^{\epsilon},-g^{\epsilon}\right)-\int_{0}^{t}\frac{e^{\sqrt{\lambda_{p}}\left(s-t\right)}}{2\sqrt{\lambda_{p}}}\left\langle f\left(s,v^{\epsilon}\left(s\right)\right),\phi_{p}\right\rangle ds\right]\phi_{p}.\label{eq:24}
\end{eqnarray}

has a unique solution $v^{\epsilon}\in C\left(\left[0,T\right];\mathcal{H}\right)$
satisfying

\begin{equation}
\left\Vert u\left(t\right)-v^{\epsilon}\left(t\right)\right\Vert \le Q\epsilon^{\frac{m\left(T-t\right)}{T}}T^{\frac{t}{T}}\left(\ln\left(\frac{T}{\epsilon^{m}}\right)\right)^{\frac{-t}{T}},\label{eq:25}
\end{equation}

where for each $p\ge1$, $\mathcal{M}_{p}:\mathcal{H}\times\mathcal{H}\to R$
such that for $w_{1},w_{2}\in\mathcal{H}$

\begin{equation}
\mathcal{M}_{p}\left(w_{1},w_{2}\right)=\left\langle w_{1},\phi_{p}\right\rangle +\frac{\left\langle w_{2},\phi_{p}\right\rangle }{\sqrt{\lambda_{p}}},\label{eq:26}
\end{equation}

and

\begin{equation}
\Phi\left(\beta,\lambda_{p},t\right)=\frac{e^{-\sqrt{\lambda_{p}}\left(T-t\right)}}{2\beta\sqrt{\lambda_{p}}+2e^{-\sqrt{\lambda_{p}}T}},\quad\Psi\left(\beta,\lambda_{p},s,t\right)=\frac{e^{-\sqrt{\lambda_{p}}\left(T+s-t\right)}}{2\beta\lambda_{p}+2\sqrt{\lambda_{p}}e^{-\sqrt{\lambda_{p}}T}},\label{eq:27}
\end{equation}

\begin{equation}
Q=\sqrt{\frac{3\lambda_{1}+3}{\lambda_{1}}}e^{\frac{3K^{2}T^{2}t}{2\lambda_{1}}}+e^{\frac{K^{2}T^{2}t}{2\lambda_{1}}}\sqrt{P}.
\end{equation}

\end{thm}
The following lemmas will lead to proof of the main theorem.
\begin{lem}
\label{lem:9}Let $\Phi\left(\beta,\lambda_{p},t\right)$ and $\Psi\left(\beta,\lambda_{p},s,t\right)$
be defined in (\ref{eq:27}), then it follows that

\begin{equation}
\Phi\left(\beta,\lambda_{p},t\right)\le\frac{1}{2}\left(\frac{\beta}{T}\right)^{\frac{-t}{T}}\left(\ln\left(\frac{T}{\beta}\right)\right)^{\frac{-t}{T}},\label{eq:29}
\end{equation}

\begin{equation}
\Psi\left(\beta,\lambda_{p},s,t\right)\le\frac{1}{2\sqrt{\lambda_{1}}}\left(\frac{\beta}{T}\right)^{\frac{s-t}{T}}\left(\ln\left(\frac{T}{\beta}\right)\right)^{\frac{s-t}{T}}.\label{eq:30}
\end{equation}

\end{lem}

\begin{lem}
\label{lem:10}The following integral equation

\begin{eqnarray}
v^{\epsilon}\left(t\right) & = & \sum_{p\ge1}\left[\Phi\left(\beta,\lambda_{p},t\right)\mathcal{M}_{p}\left(\varphi^{\epsilon},g^{\epsilon}\right)+\int_{0}^{t}\Psi\left(\beta,\lambda_{p},s,t\right)\left\langle f\left(s,v^{\epsilon}\left(s\right)\right),\phi_{p}\right\rangle ds\right]\phi_{p}\nonumber \\
 &  & +\sum_{p\ge1}\left[\frac{e^{-\sqrt{\lambda_{p}}t}}{2}\mathcal{M}_{p}\left(\varphi^{\epsilon},-g^{\epsilon}\right)-\int_{0}^{t}\frac{e^{\sqrt{\lambda_{p}}\left(s-t\right)}}{2\sqrt{\lambda_{p}}}\left\langle f\left(s,v^{\epsilon}\left(s\right)\right),\phi_{p}\right\rangle ds\right]\phi_{p},\label{eq:31}
\end{eqnarray}

has a unique solution $v^{\epsilon}\in C\left(\left[0,T\right];\mathcal{H}\right)$.
\end{lem}

\begin{lem}
\label{lem:11}The problem

\begin{eqnarray}
u^{\epsilon}\left(t\right) & = & \sum_{p\ge1}\left[\Phi\left(\beta,\lambda_{p},t\right)\mathcal{M}_{p}\left(\varphi,g\right)+\int_{0}^{t}\Psi\left(\beta,\lambda_{p},s,t\right)\left\langle f\left(s,u^{\epsilon}\left(s\right)\right),\phi_{p}\right\rangle ds\right]\phi_{p}\nonumber \\
 &  & +\sum_{p\ge1}\left[\frac{e^{-\sqrt{\lambda_{p}}t}}{2}\mathcal{M}_{p}\left(\varphi,-g\right)-\int_{0}^{t}\frac{e^{\sqrt{\lambda_{p}}\left(s-t\right)}}{2\sqrt{\lambda_{p}}}\left\langle f\left(s,u^{\epsilon}\left(s\right)\right),\phi_{p}\right\rangle ds\right]\phi_{p},\label{eq:32}
\end{eqnarray}

has a unique solution $u^{\epsilon}\in C\left(\left[0,T\right];\mathcal{H}\right)$
and the error estimate holds

\begin{equation}
\left\Vert v^{\epsilon}\left(t\right)-u^{\epsilon}\left(t\right)\right\Vert \le\sqrt{\frac{3\lambda_{1}+3}{\lambda_{1}}}e^{\frac{3K^{2}T^{2}t}{2\lambda_{1}}}\left(\frac{\beta}{T}\right)^{\frac{-t}{T}}\left(\ln\left(\frac{T}{\beta}\right)\right)^{\frac{-t}{T}}\epsilon.\label{eq:33}
\end{equation}

\end{lem}

\begin{lem}
\label{lem:12}Let $u^{\epsilon}\left(t\right)$ be a function defined
in (\ref{eq:32}), then the following estimate holds

\begin{equation}
\left\Vert u\left(t\right)-u^{\epsilon}\left(t\right)\right\Vert \le e^{\frac{T^{2}K^{2}t}{2\lambda_{1}}}\sqrt{P}\beta\left(\frac{\beta}{T}\right)^{\frac{-t}{T}}\left(\ln\left(\frac{T}{\beta}\right)\right)^{\frac{-t}{T}}.\label{eq:34}
\end{equation}

\end{lem}

\section{Numerical examples}

In this section, we aim to show two numerical examples to validate
the accuracy and efficiency of our proposed regularization method
for 1-D semi-linear elliptic problems including both linear and nonlinear
cases. The examples are involved with the operator $\mathcal{A}=-\dfrac{\partial^{2}}{\partial x^{2}}$
and taken by Hilbert space $\mathcal{H}=L^{2}\left(0,\pi\right)$.
Particularly, we give examples of a modified Helmholtz equation and
an elliptic sine-Gordon equation to demonstrate how the method works. 

The aim of numerical experiments is to observe $\epsilon=10^{-r}$
for $r\in\mathbb{N}$. The couple of $\left(\varphi^{\epsilon},g^{\epsilon}\right)$
plays as measured data with a random noise. More precisely, we take
perturbation in couple of exact data $\left(\varphi,g\right)$ to
define $\left(\varphi^{\epsilon},g^{\epsilon}\right)$ by the following
way.

\begin{eqnarray*}
\varphi^{\epsilon}\left(x\right) & = & \varphi\left(x\right)+\frac{\epsilon\cdot\mbox{rand}}{\sqrt{\pi}},\\
g^{\epsilon}\left(x\right) & = & g\left(x\right)+\frac{\epsilon\cdot\mbox{rand}}{\sqrt{\pi}},
\end{eqnarray*}

where $\mbox{rand}$ is a random number determined in $\left[-1,1\right]$.

Then, the regularized solution (with choosing $m=0.99$) is expected
to be closed to the exact solution under a proper discretization.
For convergence tests, we would like to introduce two errors: the
absolute error at the midpoint $\dfrac{\pi}{2}$ and the relative
root mean square (RRMS) error. Also, the 2-D and 3-D graphs are applied
and analysed.

To be more coherent, we are going to divide this section into two
subsections. The first one is to consider the modified Helmholtz equation
and the second one is for the elliptic sine-Gordon equation. As we
introduced, they are simply outstanding for many applied problems.
\begin{rem}
Generally, the whole process is summarized in the following steps.

\textbf{Step 1.} Given $N,K$ and $M$ to have

\[
x_{j}=j\Delta x,\Delta x=\dfrac{1}{K},j=\overline{0,K},
\]

\[
t_{i}=i\Delta t,\Delta t=\dfrac{1}{M},i=\overline{0,M}.
\]

\textbf{Step 2.} Choose $r$, put $v^{\epsilon}\left(x,t_{i}\right)=v_{i}^{\epsilon}\left(x\right),i=\overline{0,M}$
and set $v_{0}^{\epsilon}\left(x\right)=\varphi^{\epsilon}\left(x\right)$.
We find

\[
V^{\epsilon}\left(x\right)=\begin{bmatrix}v_{0}^{\epsilon}\left(x\right) & v_{1}^{\epsilon}\left(x\right) & ... & v_{M}^{\epsilon}\left(x\right)\end{bmatrix}^{T}\in\mathbb{R}^{M+1}.
\]

\textbf{Step 3.} For $i=\overline{0,M}$ and $j=\overline{0,K},$
put $v_{i}^{\epsilon}\left(x_{j}\right)=v_{i,j}^{\epsilon}$ and $u\left(x_{j},t_{i}\right)=u_{ji}$,
we find the matrices in $\mathbb{R}^{M+1}\times\mathbb{R}^{K+1}$
containing all discrete values of the exact solution $u\left(x,t\right)$
and the regularized solution $v^{\epsilon}\left(x,t\right)$, denoted
by $U$ and $V^{\epsilon}$, respectively.

\[
U=\begin{bmatrix}u_{0,0} & u_{0,1} & \cdots & u_{0,K}\\
u_{1,0} & u_{1,1} & \cdots & u_{1,K}\\
\vdots & \vdots & \ddots & \vdots\\
u_{M,0} & u_{M,1} & \ldots & u_{M,K}
\end{bmatrix},\quad V^{\epsilon}=\begin{bmatrix}v_{0,0}^{\epsilon} & v_{0,1}^{\epsilon} & \cdots & v_{0,K}^{\epsilon}\\
v_{1,0}^{\epsilon} & v_{1,1}^{\epsilon} & \cdots & v_{1,K}^{\epsilon}\\
\vdots & \vdots & \ddots & \vdots\\
v_{M,0}^{\epsilon} & v_{M,1}^{\epsilon} & \ldots & v_{M,K}^{\epsilon}
\end{bmatrix}.
\]

\textbf{Step 4.} Calculate the errors and present 2-D and 3-D graphs.

\begin{equation}
E\left(t_{i}\right)=\left|u\left(\frac{\pi}{2},t_{i}\right)-v^{\epsilon}\left(\frac{\pi}{2},t_{i}\right)\right|,\label{eq:35}
\end{equation}

\begin{equation}
R\left(t_{i}\right)=\frac{\sqrt{\sum_{0\le j\le K}\left|u\left(x_{j},t_{i}\right)-v^{\epsilon}\left(x_{j},t_{i}\right)\right|^{2}}}{\sqrt{\sum_{0\le j\le K}\left|u\left(x_{j},t_{i}\right)\right|^{2}}}.\label{eq:36}
\end{equation}

\end{rem}

\subsection{Example 1}

We will consider the following equation.

\begin{equation}
\begin{cases}
\frac{\partial^{2}}{\partial t^{2}}u\left(x,t\right)+\frac{\partial^{2}}{\partial x^{2}}u\left(x,t\right)=u\left(x,t\right) & ,\left(x,t\right)\in\left(0,\pi\right)\times\left(0,1\right),\\
\frac{\partial}{\partial x}u\left(0,t\right)=u\left(\pi,t\right)=0 & ,t\in\left(0,1\right),\\
u\left(x,0\right)=\varphi\left(x\right),\quad\frac{\partial}{\partial t}u\left(x,0\right)=0 & ,x\in\left(0,\pi\right).
\end{cases}\label{eq:37}
\end{equation}

Based on $\mathcal{D}\left(\mathcal{A}\right)=\left\{ v\in H^{1}\left(0,\pi\right):v\left(\pi\right)=0\right\} $,
we get an orthonormal eigenbasis $\phi_{p}\left(x\right)=\sqrt{\dfrac{2}{\pi}}\cos\left(\sqrt{\lambda_{p}}x\right)$
associated with the eigenvalue $\lambda_{p}=\left(p-\dfrac{1}{2}\right)^{2}$
in $L^{2}\left(0,\pi\right)$. In order to ensure the problem (\ref{eq:37})
has solution with a given Cauchy data $\varphi$, we will construct
the exact solution from a function $h$ as follows

\begin{equation}
u\left(x,1\right)=\frac{2}{\pi}\sum_{1\le p\le N}\left\langle h\left(\xi\right),\cos\left(\left(p-\frac{1}{2}\right)\xi\right)\right\rangle \cos\left(\left(p-\frac{1}{2}\right)x\right),
\end{equation}

where $N$ is a truncation term and $h$ will be chosen later. Then,
this problem has a unique solution by applying method of separation
of variables.

\begin{equation}
u\left(x,t\right)=\frac{2}{\pi}\sum_{1\le p\le N}\frac{\cosh\left(t\sqrt{\left(p-\frac{1}{2}\right)^{2}+1}\right)}{\cosh\left(\sqrt{\left(p-\frac{1}{2}\right)^{2}+1}\right)}\left\langle h\left(\xi\right),\cos\left(\left(p-\frac{1}{2}\right)\xi\right)\right\rangle \cos\left(\left(p-\frac{1}{2}\right)x\right).\label{eq:39}
\end{equation}

Thus, we have

\begin{equation}
\varphi\left(x\right)=\frac{2}{\pi}\sum_{1\le p\le N}\frac{\left\langle h\left(\xi\right),\cos\left(\left(p-\frac{1}{2}\right)\xi\right)\right\rangle }{\cosh\left(\sqrt{\left(p-\frac{1}{2}\right)^{2}+1}\right)}\cos\left(\left(p-\frac{1}{2}\right)x\right).
\end{equation}

Simultaneously, the regularized solution defined in (\ref{eq:24})
becomes

\begin{eqnarray}
v^{\epsilon}\left(x,t\right) & = & \sum_{1\le p\le N}\Phi\left(\epsilon,p,t\right)\mathcal{M}_{p}\left(\varphi^{\epsilon},g^{\epsilon}\right)\cos\left(\left(p-\frac{1}{2}\right)x\right)\nonumber \\
 &  & +\sum_{1\le p\le N}\left(\int_{0}^{t}\int_{0}^{\pi}\Psi\left(\epsilon,p,s,t\right)v^{\epsilon}\left(x,s\right)\cos\left(\left(p-\frac{1}{2}\right)x\right)dxds\right)\cos\left(\left(p-\frac{1}{2}\right)x\right)\nonumber \\
 &  & +\frac{1}{2}\sum_{1\le p\le N}e^{-\left(p-\frac{1}{2}\right)t}\mathcal{M}_{p}\left(\varphi^{\epsilon},-g^{\epsilon}\right)\cos\left(\left(p-\frac{1}{2}\right)x\right)\nonumber \\
 &  & -\sum_{1\le p\le N}\left(\frac{2}{\pi\left(2p-1\right)}\int_{0}^{t}\int_{0}^{\pi}e^{\left(p-\frac{1}{2}\right)\left(s-t\right)}v^{\epsilon}\left(x,s\right)\cos\left(\left(p-\frac{1}{2}\right)x\right)dxds\right)\cos\left(\left(p-\frac{1}{2}\right)x\right),\nonumber \\
\label{eq:41}
\end{eqnarray}

where $\mathcal{M}_{p}\left(\varphi^{\epsilon},\pm g^{\epsilon}\right),\Phi\left(\epsilon,p,t\right)$
and $\Psi\left(\epsilon,p,s,t\right)$ are induced by (\ref{eq:26})-(\ref{eq:27}).
They are explicitly defined as follows.

\begin{equation}
\mathcal{M}_{p}\left(\varphi^{\epsilon},\pm g^{\epsilon}\right)=\frac{2}{\pi}\int_{0}^{\pi}\left[\varphi^{\epsilon}\left(x\right)\pm\frac{g^{\epsilon}\left(x\right)}{p-\frac{1}{2}}\right]\cos\left(\left(p-\frac{1}{2}\right)x\right)dx,
\end{equation}

\begin{equation}
\Phi\left(\epsilon,p,t\right)=\frac{e^{-\left(p-\frac{1}{2}\right)\left(1-t\right)}}{\epsilon^{0.99}\left(2p-1\right)+2e^{-\left(p-\frac{1}{2}\right)}},\quad\Psi\left(\epsilon,p,s,t\right)=\frac{2}{\pi}\frac{e^{-\left(p-\frac{1}{2}\right)\left(1+s-t\right)}}{2\epsilon^{0.99}\left(p-\frac{1}{2}\right)^{2}+\left(2p-1\right)e^{-\left(p-\frac{1}{2}\right)}}.\label{eq:43}
\end{equation}

Now when we divide the time $t_{i}=i\Delta t,\Delta t=\dfrac{1}{M},i=\overline{0,M}$,
it turns out that a simple iterative scheme in time is applied to
(\ref{eq:41}). Particularly, we will compute $v_{i}^{\epsilon}\left(x\right),i=\overline{1,M}$
from $v_{0}^{\epsilon}\left(x\right)=\varphi^{\epsilon}\left(x\right)$
as follows.

\begin{equation}
v_{i}^{\epsilon}\left(x\right)\equiv v^{\epsilon}\left(x,t_{i}\right)=\sum_{1\le p\le N}\left[\mathcal{R}\left(\epsilon,p,t_{i}\right)-\mathcal{W}\left(\epsilon,p,t_{i}\right)\right]\cos\left(\left(p-\frac{1}{2}\right)x\right),\label{eq:44}
\end{equation}

where

\begin{eqnarray}
\mathcal{R}\left(\epsilon,p,t_{i}\right) & = & \Phi\left(\epsilon,p,t_{i}\right)\mathcal{M}_{p}\left(\varphi^{\epsilon},g^{\epsilon}\right)+\frac{1}{2}e^{-\left(p-\frac{1}{2}\right)t_{i}}\mathcal{M}_{p}\left(\varphi^{\epsilon},-g^{\epsilon}\right)\nonumber \\
 &  & +\sum_{1\le j\le i}\int_{t_{j-1}}^{t_{j}}\int_{0}^{\pi}\Psi\left(\epsilon,p,s,t_{i}\right)v_{j-1}^{\epsilon}\left(x\right)\cos\left(\left(p-\frac{1}{2}\right)x\right)dxds,\label{eq:45}
\end{eqnarray}

\begin{equation}
\mathcal{W}\left(\epsilon,p,t_{i}\right)=\frac{2}{\pi\left(2p-1\right)}\sum_{1\le j\le i}\int_{t_{j-1}}^{t_{j}}\int_{0}^{\pi}e^{\left(p-\frac{1}{2}\right)\left(s-t_{i}\right)}v_{j-1}^{\epsilon}\left(x\right)\cos\left(\left(p-\frac{1}{2}\right)x\right)dxds.\label{eq:46}
\end{equation}

As we know, $h$ plays the role as a test function. From this example,
we want to find exactly inner products between the test function and
the eigenbasis by choosing simple functions: $x^{2}\left(\pi-x\right)$
and ${\displaystyle \sum_{k=1}^{3}\frac{\cos\left(kx\right)}{k}}$.
On the other hand, we note that (\ref{eq:45}) and (\ref{eq:46})
can be simplified by directly computing the following integrations.

\begin{equation}
\int_{t_{j-1}}^{t_{j}}\Psi\left(\epsilon,p,s,t_{i}\right)ds=\frac{4}{\pi\left(1-2p\right)}\frac{e^{-\left(p-\frac{1}{2}\right)\left(1+t_{j}-t_{i}\right)}-e^{-\left(p-\frac{1}{2}\right)\left(1+t_{j-1}-t_{i}\right)}}{2\epsilon^{0.99}\left(p-\frac{1}{2}\right)^{2}+\left(2p-1\right)e^{-\left(p-\frac{1}{2}\right)}},
\end{equation}

\begin{equation}
\int_{t_{j-1}}^{t_{j}}e^{\left(p-\frac{1}{2}\right)\left(s-t_{i}\right)}ds=\frac{2}{2p-1}\left[e^{\left(p-\frac{1}{2}\right)\left(t_{j}-t_{i}\right)}-e^{\left(p-\frac{1}{2}\right)\left(t_{j-1}-t_{i}\right)}\right].
\end{equation}

\noindent \begin{center}
\begin{table}
\noindent \begin{centering}
\begin{tabular}{|c|c|c|c|c|}
\hline 
Test function & $\epsilon$ & $E\left(\dfrac{1}{10}\right)$ & $E\left(\dfrac{1}{2}\right)$ & $E\left(1\right)$\tabularnewline
\hline 
\multirow{4}{*}{$h\left(x\right)=x^{2}\left(\pi-x\right)$} & $10^{-2}$ & 0.035409705039934 & 0.116746863516900 & 0.372168953951916\tabularnewline
\cline{2-5} 
 & $10^{-4}$ & 0.000431278831272 & 0.003358920542896 & 0.023605079336301\tabularnewline
\cline{2-5} 
 & $10^{-6}$ & 0.000014949076139 & 0.001913491383348 & 0.019016995706460\tabularnewline
\cline{2-5} 
 & $10^{-8}$ & 0.000009993878017 & 0.001807682028770 & 0.018567839990525\tabularnewline
\hline 
\multirow{4}{*}{$h\left(x\right)={\displaystyle \sum_{k=1}^{3}\frac{\cos\left(kx\right)}{k}}$} & $10^{-1}$ & 0.004249941946421 & 0.074435441315929 & 0.272260206158619\tabularnewline
\cline{2-5} 
 & $10^{-3}$ & 0.001582454284463 & 0.004009212062991 & 0.014413248824993\tabularnewline
\cline{2-5} 
 & $10^{-5}$ & 0.000018911820283 & 0.000313502875558 & 0.003106592235082\tabularnewline
\cline{2-5} 
 & $10^{-7}$ & 0.000001421650997 & 0.000230738814781 & 0.002877357795009\tabularnewline
\hline 
\end{tabular}
\par\end{centering}

\caption{The absolute error at the midpoint defined in (\ref{eq:35}) with
$t=\dfrac{1}{10};\dfrac{1}{2};1$ for both two test functions in Example
1.\label{tab:1}}
\end{table}

\par\end{center}

\noindent \begin{center}
\begin{table}
\noindent \begin{centering}
\begin{tabular}{|c|c|c|c|c|c|}
\hline 
Test function & $\epsilon$ & $R\left(1\right)$ & Test function & $\epsilon$ & $R\left(1\right)$\tabularnewline
\hline 
\multirow{4}{*}{$h\left(x\right)=x^{2}\left(\pi-x\right)$} & $10^{-2}$ & 0.103782899356401 & \multirow{4}{*}{$h\left(x\right)={\displaystyle \sum_{k=1}^{3}\frac{\cos\left(kx\right)}{k}}$} & $10^{-1}$ & 0.497932025244192\tabularnewline
\cline{2-3} \cline{5-6} 
 & $10^{-4}$ & 0.005938944110216 &  & $10^{-3}$ & 0.020910786614042\tabularnewline
\cline{2-3} \cline{5-6} 
 & $10^{-6}$ & 0.004681856304455 &  & $10^{-5}$ & 0.005441953635180\tabularnewline
\cline{2-3} \cline{5-6} 
 & $10^{-8}$ & 0.004668575093985 &  & $10^{-7}$ & 0.005276479332669\tabularnewline
\hline 
\end{tabular}
\par\end{centering}

\caption{The RRMS error defined in (\ref{eq:36}) with $t=1$ for both two
test functions in Example 1.\label{tab:2}}
\end{table}

\par\end{center}

\noindent \begin{center}
\begin{figure}
\noindent \includegraphics[scale=0.6]{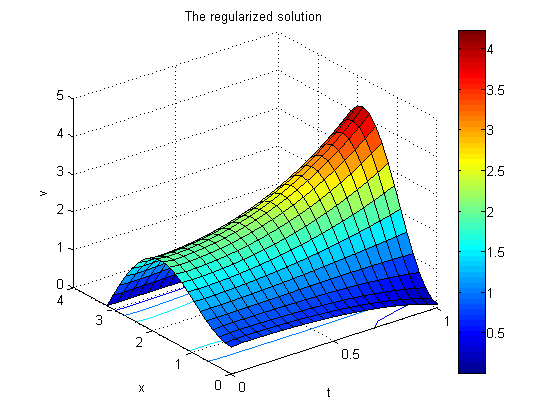}\includegraphics[scale=0.6]{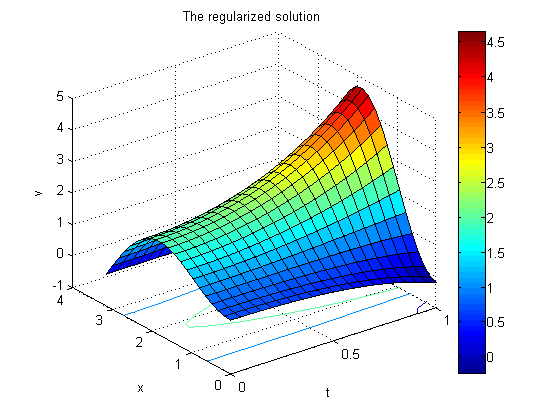}

\caption{The regularized solution (\ref{eq:41}) of Example 1 for $h\left(x\right)=x^{2}\left(\pi-x\right)$
and $\epsilon=10^{-r}$ with $r=2;4$ in 3-D representation.\label{fig:1}}
\end{figure}

\par\end{center}

\noindent \begin{center}
\begin{figure}
\noindent \includegraphics[scale=0.6]{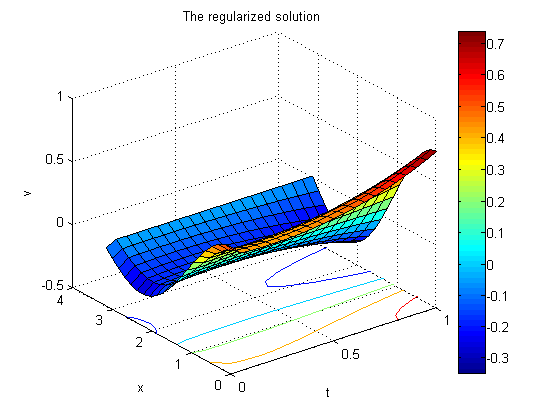}\includegraphics[scale=0.6]{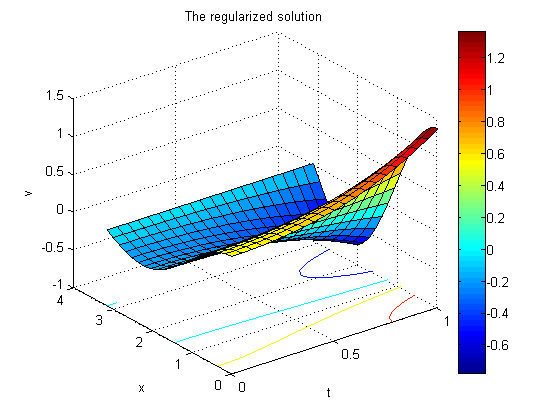}

\caption{The regularized solution (\ref{eq:41}) of Example 1 for $h\left(x\right)={\displaystyle \sum_{k=1}^{3}\frac{\cos\left(kx\right)}{k}}$
and $\epsilon=10^{-r}$ with $r=1;3$ in 3-D representation.\label{fig:2}}
\end{figure}

\par\end{center}

\begin{figure}
\noindent \includegraphics[scale=0.6]{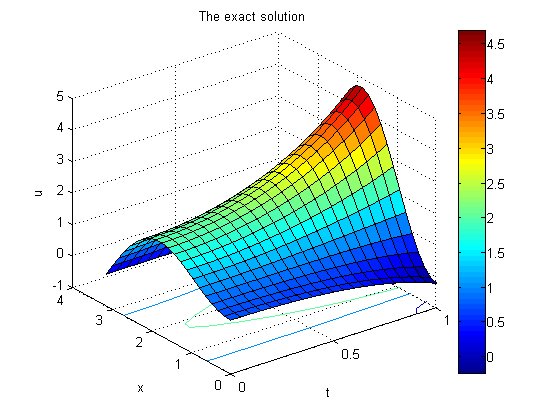}\includegraphics[scale=0.6]{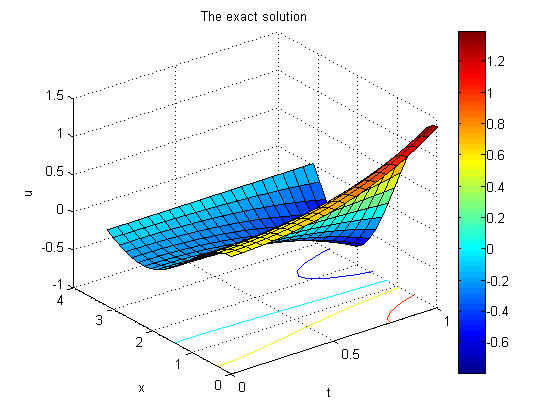}

\caption{The exact solution (\ref{eq:39}) for both two test functions in 3-D
representation in Example 1.\label{fig:3}}
\end{figure}

\begin{figure}
\noindent \includegraphics[scale=0.6]{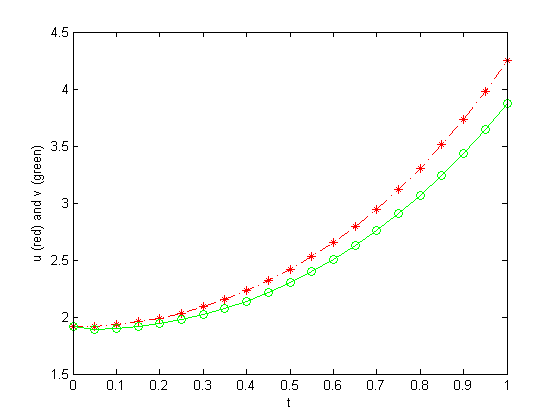}\includegraphics[scale=0.6]{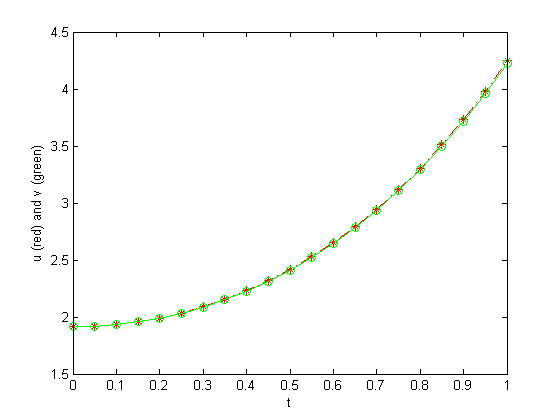}

\caption{2-D graphs of the exact solution (red) and regularized solution (green)
at $x=\dfrac{\pi}{2}$ for $h\left(x\right)=x^{2}\left(\pi-x\right)$
and $\epsilon=10^{-r}$ with $r=2;4$ in Example 1.\label{fig:4}}
\end{figure}

\begin{figure}
\noindent \includegraphics[scale=0.6]{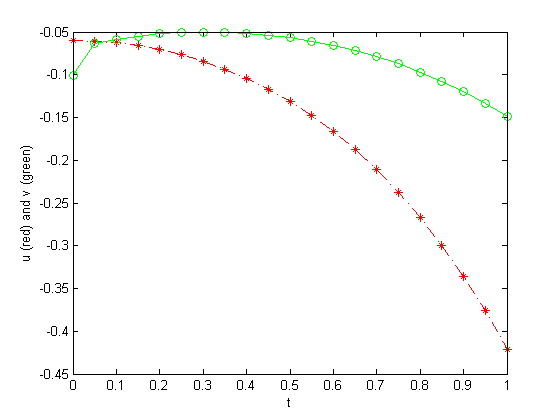}\includegraphics[scale=0.6]{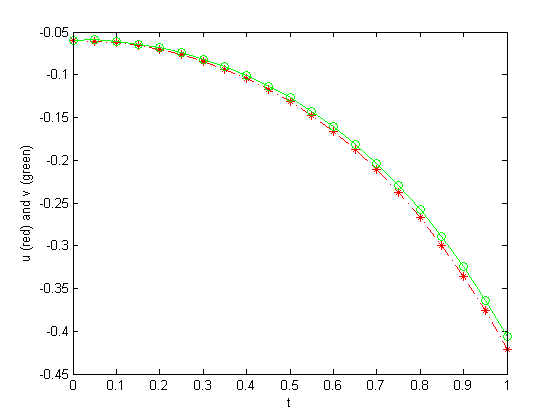}

\noindent \includegraphics[scale=0.6]{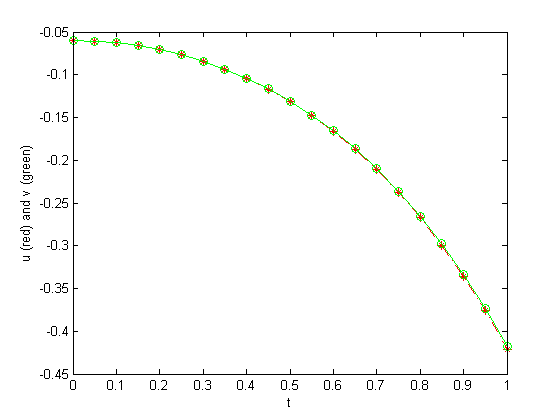}\includegraphics[scale=0.6]{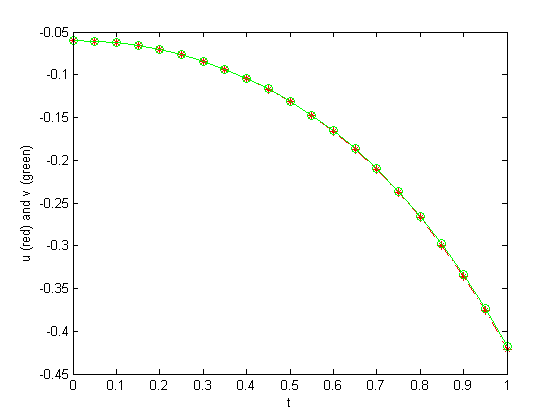}

\caption{2-D graphs of the exact solution (red) and regularized solution (green)
at $x=\dfrac{\pi}{2}$ for $h\left(x\right)={\displaystyle \sum_{k=1}^{3}\frac{\cos\left(kx\right)}{k}}$
and $\epsilon=10^{-r}$ with $r=1;3;5;7$ in Example 1.\label{fig:5}}
\end{figure}

\subsection*{Comments.}

In this computations, the square grid size for time and space variables
are rawly set by choosing $K=M=20$. The truncation term is simply
equal to $N=3$.

Table \ref{tab:1} and Table \ref{tab:2} show the absolute error
at the midpoint $\dfrac{\pi}{2}$ and RRMS error defined in (\ref{eq:35})-(\ref{eq:36})
for both two test functions $h$. Particularly, the tables show the
errors between the exact solution whose existence is ensured under
the test function $h$, recall that in this example we let $h\left(x\right)=x^{2}\left(\pi-x\right)$
and $h\left(x\right)={\displaystyle \sum_{k=1}^{3}\frac{\cos\left(kx\right)}{k}}$,
and the regularized solution (\ref{eq:41}) at the fixed time $t=\dfrac{1}{10};\dfrac{1}{2};1$
indicating three basic stage of time, nearly initial-middle-final,
are both considered. We observe that the further initial point, the
slower convergence speed and the smaller $\epsilon$, the smaller
errors. 

For the test function $h\left(x\right)=x^{2}\left(\pi-x\right)$,
we show the corresponding exact solution in Figure \ref{fig:3} (left)
and present. Despite the same 3-D shape, it should be given attention
to the color bar of the regularized ones, especially the maximum values
attaining on the bar. In addition, Figure \ref{fig:4} present the
2-D graphs of the solutions at $x=\dfrac{\pi}{2}$ for $\epsilon=10^{-2};10^{-4}$.
By observation, the regularized solution is close to the exact one
when $\epsilon$ gets smaller.

Similarly, for the test function $h\left(x\right)={\displaystyle \sum_{k=1}^{3}\frac{\cos\left(kx\right)}{k}}$
we show in Figure \ref{fig:3} (right) the exact solution and in Figure
\ref{fig:2} the regularized solution (\ref{eq:41}) for $\epsilon=10^{-1};10^{-3}$.
In Figure \ref{fig:5}, we present the 2-D graphs of the solutions
at the middle point of space for $\epsilon=10^{-r}$ with $r=1;3;5;7$,
respectively.

\subsection{Example 2}

For this example, we intend to give attention to an elliptic sine-Gordon
equation.

\begin{equation}
\begin{cases}
\frac{\partial^{2}}{\partial t^{2}}u\left(x,t\right)+\frac{\partial^{2}}{\partial x^{2}}u\left(x,t\right)=\sin\left(u\left(x,t\right)\right)-\sin\left(t\sin x\right)-t\sin x & ,\left(x,t\right)\in\left(0,\pi\right)\times\left(0,1\right),\\
u\left(0,t\right)=u\left(\pi,t\right)=0 & ,t\in\left(0,1\right),\\
u\left(x,0\right)=0,\quad\frac{\partial}{\partial t}u\left(x,0\right)=\sin x & ,x\in\left(0,\pi\right).
\end{cases}
\end{equation}

It is easy to see that for $\mathcal{D}\left(\mathcal{A}\right)=H_{0}^{1}\left(0,\pi\right)$,
we have an orthonormal eigenbasis $\phi_{p}\left(x\right)=\sqrt{\dfrac{2}{\pi}}\sin\left(\sqrt{\lambda_{p}}x\right)$
in $L^{2}\left(0,\pi\right)$ and $\lambda_{p}=p^{2}$ is the corresponding
eigenvalue. The exact solution is $u\left(x,t\right)=t\sin x$. Similar
to (\ref{eq:41})-(\ref{eq:43}) in Exampe 1, we establish the regularized
solution.

\begin{eqnarray}
v^{\epsilon}\left(x,t\right) & = & \sum_{1\le p\le N}\Phi\left(\epsilon,p,t\right)\mathcal{M}_{p}\left(\varphi^{\epsilon},g^{\epsilon}\right)\sin\left(px\right)\nonumber \\
 &  & +\sum_{1\le p\le N}\left(\int_{0}^{t}\int_{0}^{\pi}\Psi\left(\epsilon,p,s,t\right)\sin\left(v^{\epsilon}\left(x,s\right)\right)\sin\left(px\right)dxds\right)\sin\left(px\right)\nonumber \\
 &  & +\frac{1}{2}\sum_{1\le p\le N}e^{-pt}\mathcal{M}_{p}\left(\varphi^{\epsilon},-g^{\epsilon}\right)\sin\left(px\right)\nonumber \\
 &  & +\sum_{1\le p\le N}\left(\frac{1}{\pi p}\int_{0}^{t}\int_{0}^{\pi}e^{p\left(s-t\right)}\sin\left(v^{\epsilon}\left(x,s\right)\right)\sin\left(px\right)dxds\right)\sin\left(px\right),\label{eq:50}
\end{eqnarray}

where

\begin{equation}
\Phi\left(\epsilon,p,t\right)=\frac{e^{-p\left(1-t\right)}}{2p\epsilon^{0.99}+2e^{-p}},\quad\Psi\left(\epsilon,p,s,t\right)=\frac{1}{\pi}\frac{e^{-p\left(1+s-t\right)}}{p^{2}\epsilon^{0.99}+pe^{-p}},
\end{equation}

\begin{equation}
\mathcal{M}_{p}\left(\varphi^{\epsilon},\pm g^{\epsilon}\right)=\frac{2}{\pi}\int_{0}^{\pi}\left[\varphi^{\epsilon}\left(x\right)\pm\frac{g^{\epsilon}\left(x\right)}{p}\right]\sin\left(px\right)dx.
\end{equation}

In the same way, we are going to compute $v_{i}^{\epsilon}\left(x\right),i=\overline{1,M}$
from $v_{0}^{\epsilon}\left(x\right)=\varphi^{\epsilon}\left(x\right)$
as (\ref{eq:44})-(\ref{eq:46}). Consequently, the following iterative
scheme is in order.

\begin{equation}
v_{i}^{\epsilon}\left(x\right)=\sum_{1\le p\le N}\left[\mathcal{R}\left(\epsilon,p,t_{i}\right)-\mathcal{W}\left(\epsilon,p,t_{i}\right)\right]\sin\left(px\right),\label{eq:53}
\end{equation}

where

\begin{eqnarray}
\mathcal{R}\left(\epsilon,p,t_{i}\right) & = & \Phi\left(\epsilon,p,t_{i}\right)\mathcal{M}_{p}\left(\varphi^{\epsilon},g^{\epsilon}\right)+\frac{1}{2}e^{-pt_{i}}\mathcal{M}_{p}\left(\varphi^{\epsilon},-g^{\epsilon}\right)\nonumber \\
 &  & +\sum_{1\le j\le i}\int_{t_{j-1}}^{t_{j}}\int_{0}^{\pi}\Psi\left(\epsilon,p,s,t_{i}\right)\left[\sin\left(v_{j-1}^{\epsilon}\left(x\right)\right)-\sin\left(s\sin x\right)-s\sin x\right]\sin\left(px\right)dxds,\nonumber \\
\label{eq:54}
\end{eqnarray}

\begin{eqnarray}
\mathcal{W}\left(\epsilon,p,t_{i}\right) & = & \frac{1}{\pi p}\sum_{1\le j\le i}\int_{t_{j-1}}^{t_{j}}\int_{0}^{\pi}e^{p\left(s-t_{i}\right)}\left[\sin\left(v_{j-1}^{\epsilon}\left(x\right)\right)-\sin\left(s\sin x\right)-s\sin x\right]\sin\left(px\right)dxds,\label{eq:55}
\end{eqnarray}

There is a little bit marked difference in computation between (\ref{eq:54})-(\ref{eq:55})
and (\ref{eq:45})-(\ref{eq:46}). In fact, we first split $\mathcal{R}\left(\epsilon,p,t_{i}\right)$
into three appropriate terms, a term $\mathcal{R}_{1}\left(\epsilon,p,t_{i}\right)$
including $\Phi\left(\epsilon,p,t_{i}\right)\mathcal{M}_{p}\left(\varphi^{\epsilon},g^{\epsilon}\right)+\dfrac{1}{2}e^{-pt_{i}}\mathcal{M}_{p}\left(\varphi^{\epsilon},-g^{\epsilon}\right)$,
a term $\mathcal{R}_{2}\left(\epsilon,p,t_{i}\right)$ including the
nonlinearity $\sin\left(v_{j-1}^{\epsilon}\left(x\right)\right)$
and a term $\mathcal{R}_{3}\left(\epsilon,p,t_{i}\right)$ containing
the rest of this sum. In order to compute $\mathcal{R}_{2}\left(\epsilon,p,t_{i}\right)$
and $\mathcal{R}_{3}\left(\epsilon,p,t_{i}\right)$, we apply Gauss-Legendre
quadrature method (see in \cite{key-1}). In particular, we have

\begin{eqnarray}
\int_{t_{j-1}}^{t_{j}}\int_{0}^{\pi}\Psi\left(\epsilon,p,s,t_{i}\right)\sin\left(v_{j-1}^{\epsilon}\left(x\right)\right)\sin\left(px\right)dxds & = & \frac{1}{\pi}\frac{e^{-p\left(1+t_{j}-t_{i}\right)}-e^{-p\left(1+t_{j-1}-t_{i}\right)}}{p^{2}\epsilon^{0.99}+pe^{-p}}\nonumber \\
 &  & \times\sum_{r=0}^{r_{0}}\gamma_{r}\sin\left(v_{j-1}^{\epsilon}\left(x_{r}\right)\right)\sin\left(px_{r}\right),
\end{eqnarray}

\begin{eqnarray}
\int_{t_{j-1}}^{t_{j}}\int_{0}^{\pi}\Psi\left(\epsilon,p,s,t_{i}\right)\left[\sin\left(s\sin x\right)+s\sin x\right]\sin\left(px\right)dxds & = & \sum_{l=0}^{l_{0}}\sum_{r=0}^{r_{0}}\alpha_{l}\gamma_{r}\Psi\left(\epsilon,p,t_{l},t_{i}\right)\nonumber \\
 &  & \times\left[\sin\left(t_{l}\sin x_{r}\right)+t_{l}\sin x_{r}\right]\sin\left(px_{r}\right),\nonumber \\
\end{eqnarray}

where $x_{r}$ and $t_{l}$ are abscissae in $\left[0,\pi\right]$
and $\left[t_{j-1},t_{j}\right]$, respectively, and $\alpha_{l},\gamma_{r}$
are associated weights.

We also do the same way in computation of (\ref{eq:55}). Hence, (\ref{eq:53})
can be determined.

\noindent \begin{center}
\begin{table}
\noindent \begin{centering}
\begin{tabular}{|c|c|c|c|}
\hline 
$\epsilon$ & $E\left(\dfrac{1}{10}\right)$ & $E\left(\dfrac{1}{2}\right)$ & $E\left(1\right)$\tabularnewline
\hline 
$10^{-1}$ & 0.086458375926430 & 0.131568588308656 & 0.221657715167904\tabularnewline
\hline 
$10^{-2}$ & 0.005697161754899 & 0.015183097329748 & 0.056405650468800\tabularnewline
\hline 
$10^{-3}$ & 0.001067813554645 & 0.002786056926348 & 0.014399880506214\tabularnewline
\hline 
$10^{-4}$ & 0.000104838093802 & 0.000817994686682 & 0.008691680983081\tabularnewline
\hline 
$10^{-5}$ & 0.000035757102538 & 0.000617861664942 & 0.007872214352913\tabularnewline
\hline 
$10^{-6}$ & 0.000019864327276 & 0.000595051843803 & 0.007845137692661\tabularnewline
\hline 
$10^{-7}$ & 0.000017600480787 & 0.000592557387737 & 0.007844539150103\tabularnewline
\hline 
$10^{-8}$ & 0.000017498079817 & 0.000592334649247 & 0.007843831738541\tabularnewline
\hline 
\end{tabular}
\par\end{centering}

\noindent \begin{centering}
\begin{tabular}{|c|c|c|c|}
\hline 
$\epsilon$ & $R\left(\dfrac{1}{10}\right)$ & $R\left(\dfrac{1}{2}\right)$ & $R\left(1\right)$\tabularnewline
\hline 
$10^{-1}$ & 0.799075862748019 & 0.250473426345937 & 0.198371536659905\tabularnewline
\hline 
$10^{-2}$ & 0.041412415178910 & 0.021481555379548 & 0.044750235777145\tabularnewline
\hline 
$10^{-3}$ & 0.008399419222540 & 0.004819027580652 & 0.010213755307730\tabularnewline
\hline 
$10^{-4}$ & 0.000858718982259 & 0.002391935255341 & 0.006165668177321\tabularnewline
\hline 
$10^{-5}$ & 0.000627339975148 & 0.002359110404713 & 0.005644668172741\tabularnewline
\hline 
$10^{-6}$ & 0.000476784804847 & 0.002314417574336 & 0.005618874168933\tabularnewline
\hline 
$10^{-7}$ & 0.000451357224492 & 0.002306699187820 & 0.005618134790713\tabularnewline
\hline 
$10^{-8}$ & 0.000450895426251 & 0.002306548381426 & 0.005617676121443\tabularnewline
\hline 
\end{tabular}
\par\end{centering}

\caption{The absolute error at the midpoint (top) defined in (\ref{eq:35})
and RRMS error (bottom) defined in (\ref{eq:36}) at $t=\dfrac{1}{10};\dfrac{1}{2};1$
in Example 2.\label{tab:3}}
\end{table}

\par\end{center}

\begin{figure}
\includegraphics[scale=0.6]{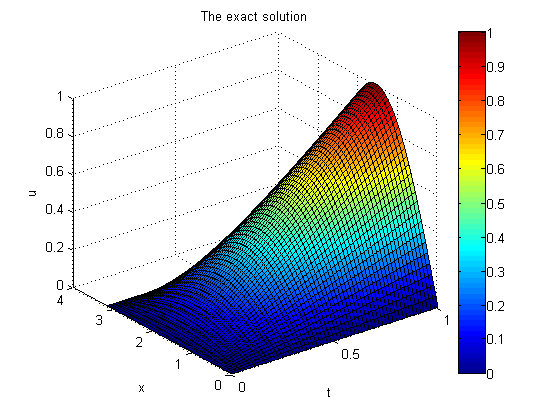}\includegraphics[scale=0.6]{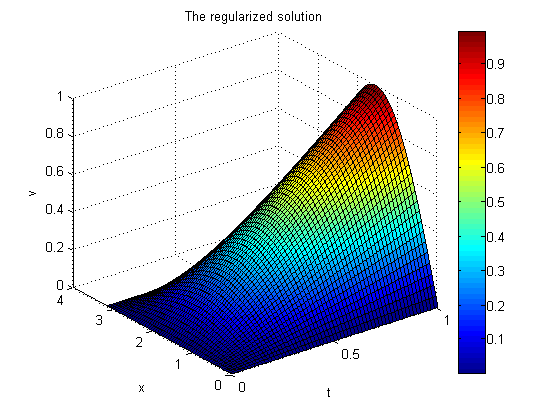}

\caption{The exact solution $u\left(x,t\right)=t\sin x$ (left) and the regularized
solution $v^{\epsilon}\left(x,t\right)$ (right) defined in (\ref{eq:50})
for $\epsilon=10^{-4}$ in 3-D representation in Example 2.\label{fig:6}}
\end{figure}

\begin{figure}
\includegraphics[scale=0.6]{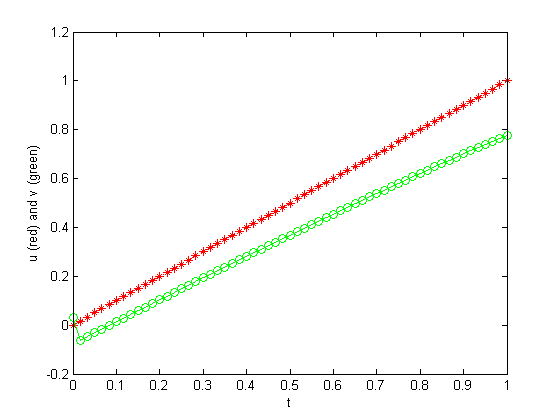}\includegraphics[scale=0.6]{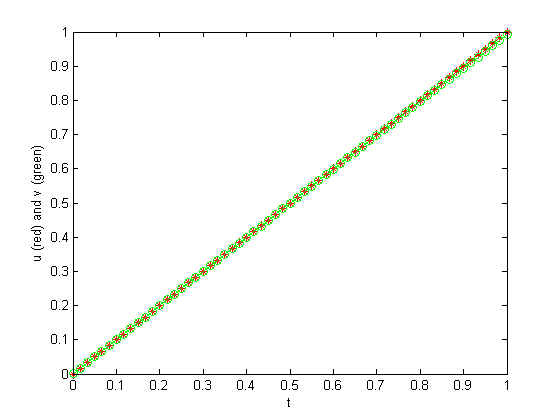}

\caption{2-D graphs of the exact solution (red) and regularized solution (green)
for $\epsilon=10^{-r}$ with $r=1;4$ in Example 2.\label{fig:7}}
\end{figure}

\subsection*{Comments.}
In this computations, the finer grid is used $\left(K=M=60\right)$
and the truncation term is still fixed as above. In the same way,
we show in Table \ref{tab:3} the errors between the exact solution
(with suppose the specific unique solution $u\left(x,t\right)=t\sin x$)
and the regularized solution (\ref{eq:50}). In Figure \ref{fig:6}
and Figure \ref{fig:7}, 3-D and 2-D graphs of them are shown, respectively.
In particular, we show in Figure \ref{fig:7} the 2-D graphs describing
how the regularized solution approachs to the exact one when $\epsilon$
becomes smaller and smaller, we illustrate the approach process by
simply presenting their graphs for $\epsilon=10^{-1};10^{-4}$. We
also show the 3-D representation of the regularized solution (with
$\epsilon=10^{-4}$) in Figure \ref{fig:6} (right).

From the numerical results, we can conclude in the same event that
the further initial point, the slower convergence speed and the convergence
is hold in general. On the other hand, it can be probably observed
that the errors reduce slowly when $\epsilon\to0$ ($\epsilon=10^{-7},10^{-8},...$),
and with a finer grid of resolution, we can have a better result in
terms of the smaller errors.

\section{Conclusion}

In this paper, we have studied the modified method to regularize the
Cauchy problem for both linear and semi-linear elliptic equations
which are severely ill-posed in general. Our approach is to present
the solution of the problem in series representation, and then propose
the regularized solution to control the strongly increasing coefficients
appearing in the series. Under some prior assumptions, we deduce error
estimates between the exact solution and regularized solution in Hilbert
space norm. The convergence rate is established by using logarithmic
estimate. We apply fundamental tools, especially using contraction
principle and Gronwall's inequality, to prove these results (see more
details in the appendix in the bottom of the paper).

In the numerical examples, we want to discuss about the semi-linear
problems with the operator $\mathcal{A}=-\Delta$ because of a wide
range of its applications. Thereby, we consider the linear Helmholtz
modified equation and the elliptic sine-Gordon equation in one-dimensional.
With lot of figures, tables and comments, our method is feasible and
efficient. The code is written in MATLAB and the computations are
done on a computer equipped with processor Pentium(R) Dual-Core CPU
2.30 GHz and having 3.0 GB total RAM.

For the other operator, the fact is that we can approximate the problem
by some numerical methods. In fact, the authors M. Charton and H.-J.
Reinhardt in \cite{key-27} apply method of lines approximation to
solve Cauchy problems for elliptic equations in two-dimensional. Particularly,
they show in this paper the approximation of $\dfrac{\partial}{\partial x}a\left(x\right)\dfrac{\partial u}{\partial x}$
under the difference schemes. Furthermore, in \cite{key-29} A. Ashyralyev
and S. Yilmaz present the first and second order of accuracy difference
schemes for the approximate solution of the initial boundary value
problem for ultra-parabolic equations with a generally positive operator.
Hence, the efficiency and feasibility of our method are obtained in
both theoretical and computational sense. It should be stated that
the issue regarding approximation of the present problem will be surveyed
in a further research.

\section*{Acknowledgments} This work is supported by Vietnam National University HoChiMinh City (VNU-HCM) under Grant No. B2014-18-01.\\
The authors would like to thank the anonymous referees for their valuable suggestions and comments leading to the improvement of our manuscript.

\section*{Appendix}

In the appendix, we would like to present the proof of all theoretical
results showed in Section 2 and Section 3 above. On account of the
proof of theorems intentionally divided into results in the related
lemmas, we will show the proof of all lemmas first, then the results
of theorems are obvious to be concluded.

\subsection*{Proof of Lemma \ref{lem:3}.}

From (\ref{eq:8})-(\ref{eq:9}), we have

\begin{eqnarray}
u^{\epsilon}\left(t\right)-v^{\epsilon}\left(t\right) & = & \sum_{p\ge1}\frac{1}{2\beta+2e^{-\sqrt{\lambda_{p}}t}}\left(\left\langle \varphi-\varphi^{\epsilon},\phi_{p}\right\rangle +\frac{\left\langle g-g^{\epsilon},\phi_{p}\right\rangle }{\sqrt{\lambda_{p}}}\right)\phi_{p}\nonumber \\
 &  & +\sum_{p\ge1}\frac{e^{-\sqrt{\lambda_{p}}t}}{2}\left(\left\langle \varphi-\varphi^{\epsilon},\phi_{p}\right\rangle -\frac{\left\langle g-g^{\epsilon},\phi_{p}\right\rangle }{\sqrt{\lambda_{p}}}\right)\phi_{p},
\end{eqnarray}

By using the inequality $\left(a+b+c+d\right)^{2}\le4\left(a^{2}+b^{2}+c^{2}+d^{2}\right)$,
we get

\begin{eqnarray}
\left|\left\langle u^{\epsilon}\left(t\right)-v^{\epsilon}\left(t\right),\phi_{p}\right\rangle \right|^{2} & \le & \frac{1}{\left(\beta+e^{-\sqrt{\lambda_{p}}t}\right)^{2}}\left(\left|\left\langle \varphi-\varphi^{\epsilon},\phi_{p}\right\rangle \right|^{2}+\frac{\left|\left\langle g-g^{\epsilon},\phi_{p}\right\rangle \right|^{2}}{\lambda_{p}}\right)\nonumber \\
 &  & +e^{-2\sqrt{\lambda_{p}}t}\left(\left|\left\langle \varphi-\varphi^{\epsilon},\phi_{p}\right\rangle \right|^{2}+\frac{\left|\left\langle g-g^{\epsilon},\phi_{p}\right\rangle \right|^{2}}{\lambda_{p}}\right).
\end{eqnarray}

Since $\beta^{2}\le\left(\beta+e^{-\sqrt{\lambda_{p}}t}\right)^{2}$
and $e^{-2\sqrt{\lambda_{p}}t}\le1\le\dfrac{1}{\beta^{2}}$, it yields 

\begin{eqnarray}
\left\Vert u^{\epsilon}\left(t\right)-v^{\epsilon}\left(t\right)\right\Vert ^{2} & = & \sum_{p\ge1}\left|\left\langle u^{\epsilon}\left(t\right)-v^{\epsilon}\left(t\right),\phi_{p}\right\rangle \right|^{2}\nonumber \\
 & \le & \frac{2}{\beta^{2}}\sum_{p\ge1}\left(\left|\left\langle \varphi-\varphi^{\epsilon},\phi_{p}\right\rangle \right|^{2}+\frac{\left|\left\langle g-g^{\epsilon},\phi_{p}\right\rangle \right|^{2}}{\lambda_{1}}\right)\nonumber \\
 & \le & \frac{2}{\beta^{2}}\left(\left\Vert \varphi-\varphi^{\epsilon}\right\Vert +\frac{\left\Vert g-g^{\epsilon}\right\Vert }{\lambda_{1}}\right).
\end{eqnarray}

Applying (\ref{eq:3}) to this, we obtain the desired result. $\square$

\subsection*{Proof of Lemma \ref{lem:4}.}

By taking the derivative of $u\left(t\right)$ in (\ref{eq:6}) with
respect to $t$, we obtain

\begin{eqnarray}
\frac{\partial}{\partial t}u\left(t\right) & = & \sum_{p\ge1}\sqrt{\lambda_{p}}\frac{e^{\sqrt{\lambda_{p}}t}}{2}\left(\left\langle \varphi-\varphi^{\epsilon},\phi_{p}\right\rangle +\frac{\left\langle g-g^{\epsilon},\phi_{p}\right\rangle }{\sqrt{\lambda_{p}}}\right)\phi_{p}\nonumber \\
 &  & -\sum_{p\ge1}\sqrt{\lambda_{p}}\frac{e^{-\sqrt{\lambda_{p}}t}}{2}\left(\left\langle \varphi-\varphi^{\epsilon},\phi_{p}\right\rangle +\frac{\left\langle g-g^{\epsilon},\phi_{p}\right\rangle }{\sqrt{\lambda_{p}}}\right)\phi_{p}.\label{eq:61}
\end{eqnarray}

It follows from (\ref{eq:6}) and (\ref{eq:61}) that

\begin{equation}
\left\langle \varphi,\phi_{p}\right\rangle +\frac{\left\langle g,\phi_{p}\right\rangle }{\sqrt{\lambda_{p}}}=e^{-\sqrt{\lambda_{p}}t}\left(\left\langle u\left(t\right),\phi_{p}\right\rangle +\frac{\left\langle \frac{\partial}{\partial t}u\left(t\right),\phi_{p}\right\rangle }{\sqrt{\lambda_{p}}}\right).
\end{equation}

Thus, we subtract $u^{\epsilon}\left(t\right)$ from $u\left(t\right)$
to have

\begin{eqnarray}
u\left(t\right)-u^{\epsilon}\left(t\right) & = & \sum_{p\ge1}\left(\frac{e^{\sqrt{\lambda_{p}}t}}{2}-\frac{1}{2\beta+2e^{-\sqrt{\lambda_{p}}t}}\right)\left(\left\langle \varphi,\phi_{p}\right\rangle +\frac{\left\langle g,\phi_{p}\right\rangle }{\sqrt{\lambda_{p}}}\right)\phi_{p}\nonumber \\
 & = & \sum_{p\ge1}\frac{\beta e^{\sqrt{\lambda_{p}}\left(t-T\right)}}{2\beta+2e^{-\sqrt{\lambda_{p}}t}}\left(\left\langle u\left(T\right),\phi_{p}\right\rangle +\frac{\left\langle \frac{\partial}{\partial t}u\left(T\right),\phi_{p}\right\rangle }{\sqrt{\lambda_{p}}}\right)\phi_{p},
\end{eqnarray}

then leads to the following

\begin{equation}
\left\Vert u\left(t\right)-u^{\epsilon}\left(t\right)\right\Vert ^{2}=\sum_{p\ge1}\left(\frac{\beta e^{\sqrt{\lambda_{p}}\left(t-T\right)}}{2\beta+2e^{-\sqrt{\lambda_{p}}t}}\right)^{2}\left(\left\langle u\left(T\right),\phi_{p}\right\rangle +\frac{\left\langle \frac{\partial}{\partial t}u\left(T\right),\phi_{p}\right\rangle }{\sqrt{\lambda_{p}}}\right)^{2}.\label{eq:64}
\end{equation}

In the next step, to get the result, we have two cases.
\begin{enumerate}
\item For $t\in\left[0,\dfrac{T}{2}\right]$, $\left\Vert u\left(t\right)-u^{\epsilon}\left(t\right)\right\Vert ^{2}$
in (\ref{eq:64}) can be estimated as follows.

\begin{eqnarray}
\left\Vert u\left(t\right)-u^{\epsilon}\left(t\right)\right\Vert ^{2} & \le & \frac{\beta^{2}}{4}\sum_{p\ge1}\left(\left\langle u\left(T\right),\phi_{p}\right\rangle +\frac{\left\langle \frac{\partial}{\partial t}u\left(T\right),\phi_{p}\right\rangle }{\sqrt{\lambda_{p}}}\right)^{2}\nonumber \\
 & \le & \frac{\beta^{2}}{2}\left(\left\Vert u\left(T\right)\right\Vert ^{2}+\frac{\left\Vert \frac{\partial}{\partial t}u\left(T\right)\right\Vert ^{2}}{\lambda_{1}}\right),
\end{eqnarray}

since $e^{\sqrt{\lambda_{p}}\left(t-T\right)}\le e^{-\sqrt{\lambda_{p}}t}$.
This implies the first estimate in (\ref{eq:17}) under condition
(\ref{eq:10}).

\item For $t\in\left[\dfrac{T}{2},T\right]$, the second estimate in (\ref{eq:17})
is obtained similarly by using the fact that $\dfrac{e^{\sqrt{\lambda_{p}}\left(t-T\right)}}{\beta+e^{-\sqrt{\lambda_{p}}t}}\le\beta^{\frac{T-t}{t}-1}$.
\end{enumerate}
Hence, we complete the proof. $\square$

\subsection*{Proof of Lemma \ref{lem:5}.}

We now rewrite the difference between $u\left(t\right)$ and $u^{\epsilon}\left(t\right)$.

\begin{equation}
u\left(t\right)-u^{\epsilon}\left(t\right)=\sum_{p\ge1}\frac{\beta}{2\beta\sqrt{\lambda_{p}}+2\sqrt{\lambda_{p}}e^{-\sqrt{\lambda_{p}}t}}\left(\sqrt{\lambda_{p}}\left\langle u\left(t\right),\phi_{p}\right\rangle +\left\langle \frac{\partial}{\partial t}u\left(t\right),\phi_{p}\right\rangle \right)\phi_{p}.\label{eq:66}
\end{equation}

We note that $e^{-\sqrt{\lambda_{p}}T}\left(\sqrt{\lambda_{p}}\left\langle u\left(T\right),\phi_{p}\right\rangle +\left\langle \frac{\partial}{\partial t}u\left(T\right),\phi_{p}\right\rangle \right)=e^{-\sqrt{\lambda_{p}}t}\left(\sqrt{\lambda_{p}}\left\langle u\left(t\right),\phi_{p}\right\rangle +\left\langle \frac{\partial}{\partial t}u\left(t\right),\phi_{p}\right\rangle \right)$,
then it follows

\begin{equation}
e^{\sqrt{\lambda_{p}}\left(T-t\right)}\left(\sqrt{\lambda_{p}}\left\langle u\left(t\right),\phi_{p}\right\rangle +\left\langle \frac{\partial}{\partial t}u\left(t\right),\phi_{p}\right\rangle \right)=\sqrt{\lambda_{p}}\left\langle u\left(T\right),\phi_{p}\right\rangle +\left\langle \frac{\partial}{\partial t}u\left(T\right),\phi_{p}\right\rangle .
\end{equation}

Thus, (\ref{eq:66}) becomes

\begin{equation}
u\left(t\right)-u^{\epsilon}\left(t\right)=\sum_{p\ge1}\frac{\beta e^{-\sqrt{\lambda_{p}}\left(T-t\right)}}{2\beta\sqrt{\lambda_{p}}+2\sqrt{\lambda_{p}}e^{-\sqrt{\lambda_{p}}t}}\left(\sqrt{\lambda_{p}}\left\langle u\left(T\right),\phi_{p}\right\rangle +\left\langle \frac{\partial}{\partial t}u\left(T\right),\phi_{p}\right\rangle \right).\label{eq:68}
\end{equation}

On the other hand, we have

\begin{equation}
\frac{\beta e^{-\sqrt{\lambda_{p}}\left(T-t\right)}}{2\beta\sqrt{\lambda_{p}}+2\sqrt{\lambda_{p}}e^{-\sqrt{\lambda_{p}}t}}\le\frac{\beta}{2\sqrt{\lambda_{1}}}\frac{e^{-\sqrt{\lambda_{p}}\left(T-t\right)}}{\frac{\beta}{\sqrt{\lambda_{1}}}\sqrt{\lambda_{p}}+e^{-\sqrt{\lambda_{p}}t}}.\label{eq:69}
\end{equation}

From (\ref{eq:68})-(\ref{eq:69}), as in proof of Lemma 4, we will
consider two cases.
\begin{enumerate}
\item For $t\in\left[0,\dfrac{T}{2}\right]$, we get

\begin{equation}
\frac{\beta e^{-\sqrt{\lambda_{p}}\left(T-t\right)}}{2\beta\sqrt{\lambda_{p}}+2\sqrt{\lambda_{p}}e^{-\sqrt{\lambda_{p}}t}}\le\frac{\beta}{2\sqrt{\lambda_{1}}}.\label{eq:70}
\end{equation}

Consequently, we obtain from (\ref{eq:68})-(\ref{eq:70})-(\ref{eq:12})
that

\begin{equation}
\left\Vert u\left(t\right)-u^{\epsilon}\left(t\right)\right\Vert ^{2}\le\frac{\beta^{2}}{4\lambda_{1}}\sum_{p\ge1}\left(\sqrt{\lambda_{p}}\left\langle u\left(T\right),\phi_{p}\right\rangle +\left\langle \frac{\partial}{\partial t}u\left(T\right),\phi_{p}\right\rangle \right)^{2}\le\frac{\beta^{2}}{4\lambda_{1}}E_{2}^{2},\label{eq:71}
\end{equation}

which implies the first estimate in (\ref{eq:18}).

\item For $t\in\left[\dfrac{T}{2},T\right]$, it follows from (\ref{eq:69})
that

\begin{eqnarray}
\frac{\beta e^{-\sqrt{\lambda_{p}}\left(T-t\right)}}{2\beta\sqrt{\lambda_{p}}+2\sqrt{\lambda_{p}}e^{-\sqrt{\lambda_{p}}t}} & \le & \frac{\beta}{2\sqrt{\lambda_{1}}}\left[\frac{T}{\frac{\beta}{\sqrt{\lambda_{1}}}\left(1+\ln\left(\frac{T}{\frac{\beta}{\sqrt{\lambda_{1}}}}\right)\right)}\right]^{\frac{2t-T}{t}}\nonumber \\
 & \le & \frac{1}{2\sqrt{\lambda_{1}}}\beta^{\frac{T-t}{t}}\left[\frac{\sqrt{\lambda_{1}}T}{1+\ln\left(\frac{\sqrt{\lambda_{1}}T}{\beta}\right)}\right]^{\frac{2t-T}{t}}.\label{eq:72}
\end{eqnarray}

Then, combining (\ref{eq:68})-(\ref{eq:72})-(\ref{eq:12}) will
give the second estimate in (\ref{eq:18}).

\end{enumerate}
Hence, we finish the proof of Lemma 5. $\square$

\subsection*{Proof of Lemma \ref{lem:6}.}

In this proof, we also obtain the estimate (\ref{eq:19}) under condition
(\ref{eq:14}) by rewriting the difference between $u\left(t\right)$
and $u^{\epsilon}\left(t\right)$,

\begin{equation}
u\left(t\right)-u^{\epsilon}\left(t\right)=\sum_{p\ge1}\frac{\beta}{2\beta+2e^{-\sqrt{\lambda_{p}}t}}\left(\left\langle u\left(t\right),\phi_{p}\right\rangle +\frac{\left\langle \frac{\partial}{\partial t}u\left(t\right),\phi_{p}\right\rangle }{\sqrt{\lambda_{p}}}\right)\phi_{p},
\end{equation}

and using a simple inequality ${\displaystyle \frac{\beta}{2\beta+2e^{-\sqrt{\lambda_{p}}t}}}\le\dfrac{\beta}{2e^{-\sqrt{\lambda_{p}}t}}$.
$\square$

\subsection*{Proof of Lemma \ref{lem:9}.}

The estimates (\ref{eq:29})-(\ref{eq:30}) are obvious under the
inequality $\dfrac{1}{\beta x+e^{-Tx}}\le\dfrac{T}{\beta\ln\left(\frac{T}{\beta}\right)}$.
Indeed, from (\ref{eq:27}) we have

\begin{eqnarray}
\Phi\left(\beta,\lambda_{p},t\right) & = & \frac{e^{-\sqrt{\lambda_{p}}\left(T-t\right)}}{2\left(\beta\sqrt{\lambda_{p}}+e^{-\sqrt{\lambda_{p}}T}\right)^{1-\frac{t}{T}}\left(\beta\sqrt{\lambda_{p}}+e^{-\sqrt{\lambda_{p}}T}\right)^{\frac{t}{T}}}\nonumber \\
 & \le & \frac{1}{2}\left(\frac{1}{\beta\sqrt{\lambda_{p}}+e^{-\sqrt{\lambda_{p}}T}}\right)^{\frac{t}{T}}\nonumber \\
 & \le & \frac{1}{2}\left(\frac{\beta}{T}\right)^{\frac{-t}{T}}\left(\ln\left(\frac{T}{\beta}\right)\right)^{\frac{-t}{T}},
\end{eqnarray}

\begin{eqnarray}
\Psi\left(\beta,\lambda_{p},s,t\right) & = & \frac{e^{-\sqrt{\lambda_{p}}\left(T-t\right)}}{2\sqrt{\lambda_{p}}\left(\beta\sqrt{\lambda_{p}}+e^{-\sqrt{\lambda_{p}}T}\right)^{1-\frac{t-s}{T}}\left(\beta\sqrt{\lambda_{p}}+e^{-\sqrt{\lambda_{p}}T}\right)^{\frac{t-s}{T}}}\nonumber \\
 & \le & \frac{1}{2}\left(\frac{1}{\beta\sqrt{\lambda_{p}}+e^{-\sqrt{\lambda_{p}}T}}\right)^{\frac{t-s}{T}}\nonumber \\
 & \le & \frac{1}{2\sqrt{\lambda_{1}}}\left(\frac{\beta}{T}\right)^{\frac{s-t}{T}}\left(\ln\left(\frac{T}{\beta}\right)\right)^{\frac{s-t}{T}}.
\end{eqnarray}

Therefore, the proof is completed. $\square$

\subsection*{Proof of Lemma \ref{lem:10}.}

For $w\in C\left(\left[0,T\right];\mathcal{H}\right)$, we consider
the following function

\begin{eqnarray}
F\left(w\right)\left(t\right) & = & \sum_{p\ge1}\left[\Phi\left(\beta,\lambda_{p},t\right)\mathcal{M}_{p}\left(w_{1},w_{2}\right)+\int_{0}^{t}\Psi\left(\beta,\lambda_{p},s,t\right)\left\langle f\left(s,w\left(s\right)\right),\phi_{p}\right\rangle ds\right]\phi_{p}\nonumber \\
 &  & +\sum_{p\ge1}\left[\frac{e^{-\sqrt{\lambda_{p}}t}}{2}\mathcal{M}_{p}\left(w_{1},-w_{2}\right)-\int_{0}^{t}\frac{e^{\sqrt{\lambda_{p}}\left(s-t\right)}}{2\sqrt{\lambda_{p}}}\left\langle f\left(s,w\left(s\right)\right),\phi_{p}\right\rangle ds\right]\phi_{p}.
\end{eqnarray}

By defining 

\begin{equation}
\Lambda\equiv\Lambda\left(\beta,\lambda_{p},t,w_{1},w_{2}\right)=\Phi\left(\beta,\lambda_{p},t\right)\mathcal{M}_{p}\left(w_{1},w_{2}\right)+\frac{e^{-\sqrt{\lambda_{p}}t}}{2}\mathcal{M}_{p}\left(w_{1},-w_{2}\right),
\end{equation}

$F\left(w\right)\left(t\right)$ becomes

\begin{equation}
F\left(w\right)\left(t\right)=\sum_{p\ge1}\left[\Lambda+\int_{0}^{t}\left(\Psi\left(\beta,\lambda_{p},s,t\right)-\frac{e^{\sqrt{\lambda_{p}}\left(s-t\right)}}{2\sqrt{\lambda_{p}}}\right)\left\langle f\left(s,w\left(s\right)\right),\phi_{p}\right\rangle ds\right]\phi_{p}.
\end{equation}

We claim that, for every $w,v\in C\left(\left[0,T\right];\mathcal{H}\right)$
and $m\ge1$, we have

\begin{equation}
\left\Vert F^{m}\left(w\right)\left(t\right)-F^{m}\left(v\right)\left(t\right)\right\Vert ^{2}\le\left(\frac{T^{3}K^{2}\beta^{-2}}{\lambda_{1}}\right)^{m}\frac{t^{m}}{m!}\left|\left\Vert w-v\right\Vert \right|^{2},\label{eq:79}
\end{equation}

where $\left|\left\Vert .\right\Vert \right|$ is supremum norm in
$C\left(\left[0,T\right];\mathcal{H}\right)$. We shall prove this
inequality by induction. Indeed, for $m=1$, we get the following
estimate.

\begin{eqnarray}
\left\Vert F^{m}\left(w\right)\left(t\right)-F^{m}\left(v\right)\left(t\right)\right\Vert ^{2} & = & \sum_{p\ge1}\left[\int_{0}^{t}\left(\Psi\left(\beta,\lambda_{p},s,t\right)-\frac{e^{\sqrt{\lambda_{p}}\left(s-t\right)}}{2\sqrt{\lambda_{p}}}\right)\left\langle f\left(s,w\left(s\right)\right)-f\left(s,v\left(s\right)\right),\phi_{p}\right\rangle ds\right]^{2}\nonumber \\
 & \le & \sum_{p\ge1}\int_{0}^{t}\left(\Psi\left(\beta,\lambda_{p},s,t\right)-\frac{e^{\sqrt{\lambda_{p}}\left(s-t\right)}}{2\sqrt{\lambda_{p}}}\right)^{2}ds\int_{0}^{t}\left|\left\langle f\left(s,w\left(s\right)\right)-f\left(s,v\left(s\right)\right),\phi_{p}\right\rangle \right|^{2}ds.\nonumber \\
\end{eqnarray}

Using the following estimate

\begin{eqnarray}
\left(\Psi\left(\beta,\lambda_{p},s,t\right)-\frac{e^{\sqrt{\lambda_{p}}\left(s-t\right)}}{2\sqrt{\lambda_{p}}}\right)^{2} & \le & 2\Psi^{2}\left(\beta,\lambda_{p},s,t\right)+\frac{e^{2\sqrt{\lambda_{p}}\left(s-t\right)}}{2\sqrt{\lambda_{p}}}\nonumber \\
 & \le & 2\Psi^{2}\left(\beta,\lambda_{p},s,t\right)+\frac{1}{2\lambda_{1}}\nonumber \\
 & \le & 2\frac{1}{4\lambda_{1}}\left(\frac{\beta}{T}\right)^{\frac{2s-2t}{T}}\left(\ln\left(\frac{T}{\beta}\right)\right)^{\frac{2s-2t}{T}}+\frac{1}{2\lambda_{1}}\nonumber \\
 & \le & \frac{1}{\lambda_{1}}\left(\frac{\beta}{T}\right)^{-2},
\end{eqnarray}

we thus have

\begin{eqnarray}
\left\Vert F^{m}\left(w\right)\left(t\right)-F^{m}\left(v\right)\left(t\right)\right\Vert ^{2} & \le & \frac{1}{\lambda_{1}}T^{2}\beta^{-2}t\int_{0}^{t}\left\Vert f\left(s,w\left(s\right)\right)-f\left(s,v\left(s\right)\right)\right\Vert ^{2}ds\nonumber \\
 & \le & \frac{1}{\lambda_{1}}T^{2}\beta^{-2}K^{2}t\int_{0}^{t}\left\Vert w\left(s\right)-v\left(s\right)\right\Vert ^{2}ds\nonumber \\
 & \le & \frac{T^{3}K^{2}\beta^{-2}}{\lambda_{1}}t\left|\left\Vert w-v\right\Vert \right|^{2}.
\end{eqnarray}

Thus (\ref{eq:79}) holds for $m=1$. Next, suppose that (\ref{eq:79})
holds for $m=k$, we prove that (\ref{eq:79}) also holds for $m=k+1$.
We have

\begin{eqnarray}
\left\Vert F^{k+1}\left(w\right)\left(t\right)-F^{k+1}\left(v\right)\left(t\right)\right\Vert ^{2} & \le & \frac{1}{\lambda_{1}}T^{2}\beta^{-2}t\int_{0}^{t}\left\Vert f\left(s,F^{k}\left(w\right)\left(s\right)\right)-f\left(s,F^{k}\left(v\right)\left(s\right)\right)\right\Vert ^{2}ds\nonumber \\
 & \le & \frac{1}{\lambda_{1}}T^{3}\beta^{-2}K^{2}\int_{0}^{t}\left(\frac{T^{3}K^{2}\beta^{-2}}{\lambda_{1}}\right)^{k}\frac{s^{k}}{k!}\left|\left\Vert w-v\right\Vert \right|^{2}ds\nonumber \\
 & \le & \left(\frac{T^{3}K^{2}\beta^{-2}}{\lambda_{1}}\right)^{k+1}\frac{t^{k+1}}{\left(k+1\right)!}\left|\left\Vert w-v\right\Vert \right|^{2}.
\end{eqnarray}

Therefore, by the induction principle, we obtain

\begin{equation}
\left\Vert F^{m}\left(w\right)\left(t\right)-F^{m}\left(v\right)\left(t\right)\right\Vert \le\sqrt{\left(\frac{T^{3}K^{2}\beta^{-2}}{\lambda_{1}}\right)^{m}\frac{t^{m}}{m!}}\left|\left\Vert w-v\right\Vert \right|,\label{eq:84}
\end{equation}

for all $w,v\in C\left(\left[0,T\right];\mathcal{H}\right)$.

We consider $F:C\left(\left[0,T\right];\mathcal{H}\right)\to C\left(\left[0,T\right];\mathcal{H}\right)$
and may see that

\[
\lim_{m\to\infty}\sqrt{\left(\frac{T^{3}K^{2}\beta^{-2}}{\lambda_{1}}\right)^{m}\frac{t^{m}}{m!}}=0.
\]

Thus, there exists a positive integer number $m_{0}$ such that 

\[
\sqrt{\left(\frac{T^{3}K^{2}\beta^{-2}}{\lambda_{1}}\right)^{m_{0}}\frac{t^{m_{0}}}{m_{0}!}}<1,
\]

and $F^{m_{0}}$ is a contraction indicating the equation $F^{m_{0}}\left(w\right)=w$
has a unique solution $w\in C\left(\left[0,T\right];\mathcal{H}\right)$.
Moreover, the fact is that $F\left(F^{m_{0}}\left(w\right)\right)=F\left(w\right)$,
then $F^{m_{0}}\left(F\left(w\right)\right)=F\left(w\right)$. By
the uniqueness of the fixed point of $F^{m_{0}}$, the equation $F\left(w\right)=w$
has a unique solution in $C\left(\left[0,T\right];\mathcal{H}\right)$. 

Hence, we obtain the result of this lemma. $\square$

\subsection*{Proof of Lemma \ref{lem:11}.}

From (\ref{eq:24}) and (\ref{eq:32}), it is clear that

\begin{eqnarray}
v^{\epsilon}\left(t\right)-u^{\epsilon}\left(t\right) & = & \sum_{p\ge1}\left[\Phi\left(\beta,\lambda_{p},t\right)\mathcal{M}_{p}\left(\varphi^{\epsilon}-\varphi,g^{\epsilon}-g\right)+\int_{0}^{t}\Psi\left(\beta,\lambda_{p},s,t\right)\left\langle f\left(s,v^{\epsilon}\left(s\right)\right)-f\left(s,u^{\epsilon}\left(s\right)\right),\phi_{p}\right\rangle ds\right]\phi_{p}\nonumber \\
 &  & +\sum_{p\ge1}\left[\frac{e^{-\sqrt{\lambda_{p}}t}}{2}\mathcal{M}_{p}\left(\varphi^{\epsilon}-\varphi,g-g^{\epsilon}\right)-\int_{0}^{t}\frac{e^{\sqrt{\lambda_{p}}\left(s-t\right)}}{2\sqrt{\lambda_{p}}}\left\langle f\left(s,v^{\epsilon}\left(s\right)\right)-f\left(s,u^{\epsilon}\left(s\right)\right),\phi_{p}\right\rangle ds\right]\phi_{p}.\nonumber \\
\label{eq:85}
\end{eqnarray}

Now we put

\begin{equation}
\eta_{1}\equiv\eta_{1}\left(\beta,\lambda_{p},t,\epsilon\right)=\Phi\left(\beta,\lambda_{p},t\right)\mathcal{M}_{p}\left(\varphi^{\epsilon}-\varphi,g^{\epsilon}-g\right),
\end{equation}

\begin{equation}
\eta_{2}\equiv\eta_{2}\left(\lambda_{p},t,\epsilon\right)=\frac{e^{-\sqrt{\lambda_{p}}t}}{2}\mathcal{M}_{p}\left(\varphi^{\epsilon}-\varphi,g-g^{\epsilon}\right),
\end{equation}

\begin{equation}
\eta_{3}\equiv\eta_{3}\left(\beta,\lambda_{p},s,t,\epsilon\right)=\int_{0}^{t}\left(\Psi\left(\beta,\lambda_{p},s,t\right)-\frac{e^{\sqrt{\lambda_{p}}\left(s-t\right)}}{2\sqrt{\lambda_{p}}}\right)\left\langle f\left(s,v^{\epsilon}\left(s\right)\right)-f\left(s,u^{\epsilon}\left(s\right)\right),\phi_{p}\right\rangle ds.
\end{equation}

We shall estimate these terms as follows. First, by (\ref{eq:26})
and (\ref{eq:29}) $\eta_{1}$ can be estimated in the following way.

\begin{eqnarray}
\eta_{1}^{2} & \le & \frac{1}{4}\left(\frac{\beta}{T}\right)^{\frac{-2t}{T}}\left(\ln\left(\frac{T}{\beta}\right)\right)^{\frac{-2t}{T}}\left(\left\langle \varphi^{\epsilon}-\varphi,\phi_{p}\right\rangle +\frac{\left\langle g^{\epsilon}-g,\phi_{p}\right\rangle }{\sqrt{\lambda_{p}}}\right)^{2}\nonumber \\
 & \le & \frac{1}{2}\left(\frac{\beta}{T}\right)^{\frac{-2t}{T}}\left(\ln\left(\frac{T}{\beta}\right)\right)^{\frac{-2t}{T}}\left(\left|\left\langle \varphi^{\epsilon}-\varphi,\phi_{p}\right\rangle \right|^{2}+\frac{\left|\left\langle g^{\epsilon}-g,\phi_{p}\right\rangle \right|^{2}}{\lambda_{1}}\right).\label{eq:89}
\end{eqnarray}

Second, we apply (\ref{eq:26}) and use the inequality $\left(\dfrac{\beta}{T}\right)^{\frac{-2t}{T}}\left(\ln\left(\dfrac{T}{\beta}\right)\right)^{\frac{-2t}{T}}\ge1$
to obtain the estimate of $\eta_{2}$.

\begin{equation}
\eta_{2}^{2}\le\frac{1}{2}\left(\frac{\beta}{T}\right)^{\frac{-2t}{T}}\left(\ln\left(\frac{T}{\beta}\right)\right)^{\frac{-2t}{T}}\left(\left|\left\langle \varphi^{\epsilon}-\varphi,\phi_{p}\right\rangle \right|^{2}+\frac{\left|\left\langle g^{\epsilon}-g,\phi_{p}\right\rangle \right|^{2}}{\lambda_{1}}\right).\label{eq:90}
\end{equation}

Finally, since (\ref{eq:30}), we get the estimate of $\eta_{3}$.

\begin{eqnarray}
\eta_{3}^{2} & \le & t^{2}\int_{0}^{t}\left(\Psi\left(\beta,\lambda_{p},s,t\right)-\frac{e^{\sqrt{\lambda_{p}}\left(s-t\right)}}{2\sqrt{\lambda_{p}}}\right)^{2}\left|\left\langle f\left(s,v^{\epsilon}\left(s\right)\right)-f\left(s,u^{\epsilon}\left(s\right)\right),\phi_{p}\right\rangle \right|^{2}ds\nonumber \\
 & \le & T^{2}\int_{0}^{t}\left(2\Psi^{2}\left(\beta,\lambda_{p},s,t\right)+\frac{e^{2\sqrt{\lambda_{p}}\left(s-t\right)}}{2\lambda_{p}}\right)\left|\left\langle f\left(s,v^{\epsilon}\left(s\right)\right)-f\left(s,u^{\epsilon}\left(s\right)\right),\phi_{p}\right\rangle \right|^{2}ds\nonumber \\
 & \le & T^{2}\int_{0}^{t}\left(\frac{1}{2\lambda_{1}}\left(\frac{\beta}{T}\right)^{\frac{2s-2t}{T}}\left(\ln\left(\frac{T}{\beta}\right)\right)^{\frac{2s-2t}{T}}+\frac{1}{2\lambda_{1}}\right)\left|\left\langle f\left(s,v^{\epsilon}\left(s\right)\right)-f\left(s,u^{\epsilon}\left(s\right)\right),\phi_{p}\right\rangle \right|^{2}ds\nonumber \\
 & \le & \frac{T^{2}}{\lambda_{1}}\int_{0}^{t}\left(\frac{\beta}{T}\right)^{\frac{2s-2t}{T}}\left(\ln\left(\frac{T}{\beta}\right)\right)^{\frac{2s-2t}{T}}\left|\left\langle f\left(s,v^{\epsilon}\left(s\right)\right)-f\left(s,u^{\epsilon}\left(s\right)\right),\phi_{p}\right\rangle \right|^{2}ds.\label{eq:91}
\end{eqnarray}

It follows from (\ref{eq:85}) and (\ref{eq:89})-(\ref{eq:91}) that

\begin{eqnarray}
\left\Vert v^{\epsilon}\left(t\right)-u^{\epsilon}\left(t\right)\right\Vert ^{2} & \le & 3\sum_{p\ge1}\left(\eta_{1}^{2}+\eta_{2}^{2}+\eta_{3}^{2}\right)\nonumber \\
 & \le & 3\left(\frac{\beta}{T}\right)^{\frac{-2t}{T}}\left(\ln\left(\frac{T}{\beta}\right)\right)^{\frac{-2t}{T}}\sum_{p\ge1}\left(\left|\left\langle \varphi^{\epsilon}-\varphi,\phi_{p}\right\rangle \right|^{2}+\frac{\left|\left\langle g^{\epsilon}-g,\phi_{p}\right\rangle \right|^{2}}{\lambda_{1}}\right)\nonumber \\
 &  & +\frac{3T^{2}}{\lambda_{1}}\sum_{p\ge1}\int_{0}^{t}\left(\frac{\beta}{T}\right)^{\frac{2s-2t}{T}}\left(\ln\left(\frac{T}{\beta}\right)\right)^{\frac{2s-2t}{T}}\left|\left\langle f\left(s,v^{\epsilon}\left(s\right)\right)-f\left(s,u^{\epsilon}\left(s\right)\right),\phi_{p}\right\rangle \right|^{2}ds.\label{eq:92}
\end{eqnarray}

Because of the fact that

\begin{eqnarray}
\sum_{p\ge1}\left(\left|\left\langle \varphi^{\epsilon}-\varphi,\phi_{p}\right\rangle \right|^{2}+\frac{\left|\left\langle g^{\epsilon}-g,\phi_{p}\right\rangle \right|^{2}}{\lambda_{1}}\right) & = & \left\Vert \varphi^{\epsilon}-\varphi\right\Vert ^{2}+\frac{\left\Vert g^{\epsilon}-g\right\Vert ^{2}}{\lambda_{1}}\nonumber \\
 & \le & \left(1+\frac{1}{\lambda_{1}}\right)\epsilon^{2},
\end{eqnarray}

we continue to get from (\ref{eq:92}) that

\begin{eqnarray}
\left\Vert v^{\epsilon}\left(t\right)-u^{\epsilon}\left(t\right)\right\Vert ^{2} & \le & 3\left(1+\frac{1}{\lambda_{1}}\right)\left(\frac{\beta}{T}\right)^{\frac{-2t}{T}}\left(\ln\left(\frac{T}{\beta}\right)\right)^{\frac{-2t}{T}}\epsilon^{2}\nonumber \\
 &  & +\frac{3T^{2}}{\lambda_{1}}\int_{0}^{t}\left(\frac{\beta}{T}\right)^{\frac{2s-2t}{T}}\left(\ln\left(\frac{T}{\beta}\right)\right)^{\frac{2s-2t}{T}}\left\Vert f\left(s,v^{\epsilon}\left(s\right)\right)-f\left(s,u^{\epsilon}\left(s\right)\right)\right\Vert ^{2}ds\nonumber \\
 & \le & 3\left(1+\frac{1}{\lambda_{1}}\right)\left(\frac{\beta}{T}\right)^{\frac{-2t}{T}}\left(\ln\left(\frac{T}{\beta}\right)\right)^{\frac{-2t}{T}}\epsilon^{2}\nonumber \\
 &  & +\frac{3K^{2}T^{2}}{\lambda_{1}}\int_{0}^{t}\left(\frac{\beta}{T}\right)^{\frac{2s-2t}{T}}\left(\ln\left(\frac{T}{\beta}\right)\right)^{\frac{2s-2t}{T}}\left\Vert v^{\epsilon}\left(s\right)-u^{\epsilon}\left(s\right)\right\Vert ^{2}ds.\label{eq:94}
\end{eqnarray}

Multiplying both sides of (\ref{eq:94}) by ${\displaystyle \left(\frac{\beta}{T}\right)^{\frac{2t}{T}}\left(\ln\left(\frac{T}{\beta}\right)\right)^{\frac{2t}{T}}}$,
it yields

\begin{eqnarray*}
\left(\frac{\beta}{T}\right)^{\frac{2t}{T}}\left(\ln\left(\frac{T}{\beta}\right)\right)^{\frac{2t}{T}}\left\Vert v^{\epsilon}\left(t\right)-u^{\epsilon}\left(t\right)\right\Vert ^{2} & \le & 3\left(1+\frac{1}{\lambda_{1}}\right)\epsilon^{2}\\
 &  & +\frac{3K^{2}T^{2}}{\lambda_{1}}\int_{0}^{t}\left(\frac{\beta}{T}\right)^{\frac{2s}{T}}\left(\ln\left(\frac{T}{\beta}\right)\right)^{\frac{2s}{T}}\left\Vert v^{\epsilon}\left(s\right)-u^{\epsilon}\left(s\right)\right\Vert ^{2}ds.
\end{eqnarray*}

By using Gronwall's inequality, we thus obtain

\begin{equation}
\left(\frac{\beta}{T}\right)^{\frac{2t}{T}}\left(\ln\left(\frac{T}{\beta}\right)\right)^{\frac{2t}{T}}\left\Vert v^{\epsilon}\left(t\right)-u^{\epsilon}\left(t\right)\right\Vert ^{2}\le3e^{\frac{3K^{2}T^{2}t}{\lambda_{1}}}\left(1+\frac{1}{\lambda_{1}}\right)\epsilon^{2},
\end{equation}

which gives the desired result (\ref{eq:33}). $\square$

\subsection*{Proof of Lemma \ref{lem:12}.}

By taking the derivative of $u\left(t\right)$ in (\ref{eq:7}) with
respect to $t$, we have

\begin{eqnarray}
\frac{\partial}{\partial t}u\left(t\right) & = & \sum_{p\ge1}\sqrt{\lambda_{p}}\left[\frac{e^{\sqrt{\lambda_{p}}t}}{2}\left(\left\langle \varphi,\phi_{p}\right\rangle +\frac{\left\langle g,\phi_{p}\right\rangle }{\sqrt{\lambda_{p}}}\right)+\int_{0}^{t}\frac{e^{\sqrt{\lambda_{p}}\left(t-s\right)}}{2\sqrt{\lambda_{p}}}\left\langle f\left(s,u\left(s\right)\right),\phi_{p}\right\rangle ds\right]\phi_{p}\nonumber \\
 &  & -\sum_{p\ge1}\sqrt{\lambda_{p}}\left[\frac{e^{\sqrt{\lambda_{p}}t}}{2}\left(\left\langle \varphi,\phi_{p}\right\rangle -\frac{\left\langle g,\phi_{p}\right\rangle }{\sqrt{\lambda_{p}}}\right)-\int_{0}^{t}\frac{e^{\sqrt{\lambda_{p}}\left(s-t\right)}}{2\sqrt{\lambda_{p}}}\left\langle f\left(s,u\left(s\right)\right),\phi_{p}\right\rangle ds\right]\phi_{p}.\label{eq:96}
\end{eqnarray}

This follows that

\begin{eqnarray}
\left\langle u\left(t\right),\phi_{p}\right\rangle +\frac{\left\langle \frac{\partial}{\partial t}u\left(t\right),\phi_{p}\right\rangle }{\sqrt{\lambda_{p}}} & = & e^{\sqrt{\lambda_{p}}t}\left(\left\langle \varphi,\phi_{p}\right\rangle +\frac{\left\langle g,\phi_{p}\right\rangle }{\sqrt{\lambda_{p}}}+\int_{0}^{t}\frac{e^{-\sqrt{\lambda_{p}}s}}{\sqrt{\lambda_{p}}}\left\langle f\left(s,u\left(s\right)\right),\phi_{p}\right\rangle ds\right)\nonumber \\
 & = & \frac{e^{\sqrt{\lambda_{p}}t}}{\sqrt{\lambda_{p}}}\left(\sqrt{\lambda_{p}}\mathcal{M}\left(\varphi,g\right)+\int_{0}^{t}e^{-\sqrt{\lambda_{p}}s}\left\langle f\left(s,u\left(s\right)\right),\phi_{p}\right\rangle ds\right).\label{eq:97}
\end{eqnarray}

Let us return to the formula of $u^{\epsilon}\left(t\right)$ in (\ref{eq:32}),
then subtracting $u^{\epsilon}\left(t\right)$ from $u\left(t\right)$,
using (\ref{eq:97}) and having direct computation yield

\begin{eqnarray}
u\left(t\right)-u^{\epsilon}\left(t\right) & = & \sum_{p\ge1}\frac{\beta\sqrt{\lambda_{p}}}{2\beta\sqrt{\lambda_{p}}+2e^{-\sqrt{\lambda_{p}}T}}\left(\left\langle u\left(t\right),\phi_{p}\right\rangle +\frac{\left\langle \frac{\partial}{\partial t}u\left(t\right),\phi_{p}\right\rangle }{\sqrt{\lambda_{p}}}\right)\phi_{p}\nonumber \\
 &  & +\sum_{p\ge1}\left[\int_{0}^{t}\left(\Psi\left(\beta,\lambda_{p},s,t\right)-\frac{e^{\sqrt{\lambda_{p}}\left(s-t\right)}}{2\sqrt{\lambda_{p}}}\right)\left\langle f\left(s,u\left(s\right)\right)-f\left(s,u^{\epsilon}\left(s\right)\right),\phi_{p}\right\rangle ds\right]\phi_{p}.
\end{eqnarray}

We thus have

\begin{eqnarray}
\left\Vert u\left(t\right)-u^{\epsilon}\left(t\right)\right\Vert ^{2} & = & \sum_{p\ge1}\left(\frac{\beta}{2\beta\sqrt{\lambda_{p}}+2e^{-\sqrt{\lambda_{p}}T}}\right)^{2}\left(\sqrt{\lambda_{p}}\left\langle u\left(t\right),\phi_{p}\right\rangle +\left\langle \frac{\partial}{\partial t}u\left(t\right),\phi_{p}\right\rangle \right)^{2}\nonumber \\
 &  & +\sum_{p\ge1}\left[\int_{0}^{t}\left(\Psi\left(\beta,\lambda_{p},s,t\right)-\frac{e^{\sqrt{\lambda_{p}}\left(s-t\right)}}{2\sqrt{\lambda_{p}}}\right)\left\langle f\left(s,u\left(s\right)\right)-f\left(s,u^{\epsilon}\left(s\right)\right),\phi_{p}\right\rangle ds\right]^{2}\nonumber \\
 & \le & \beta^{2}\sum_{p\ge1}\Phi^{2}\left(\epsilon,\lambda_{p},t\right)e^{2\sqrt{\lambda_{p}}\left(T-t\right)}\left(\sqrt{\lambda_{p}}\left\langle u\left(t\right),\phi_{p}\right\rangle +\left\langle \frac{\partial}{\partial t}u\left(t\right),\phi_{p}\right\rangle \right)^{2}\nonumber \\
 &  & +T^{2}\sum_{p\ge1}\int_{0}^{t}\left(\Psi\left(\beta,\lambda_{p},s,t\right)-\frac{e^{\sqrt{\lambda_{p}}\left(s-t\right)}}{2\sqrt{\lambda_{p}}}\right)^{2}\left|\left\langle f\left(s,u\left(s\right)\right)-f\left(s,u^{\epsilon}\left(s\right)\right),\phi_{p}\right\rangle \right|^{2}ds.\label{eq:99}
\end{eqnarray}

Now we put $\rho_{1},\rho_{2}$ as

\begin{equation}
\rho_{1}\left(\beta,\lambda_{p},t,\epsilon\right)=\beta^{2}\sum_{p\ge1}\Phi^{2}\left(\epsilon,\lambda_{p},t\right)e^{2\sqrt{\lambda_{p}}\left(T-t\right)}\left(\sqrt{\lambda_{p}}\left\langle u\left(t\right),\phi_{p}\right\rangle +\left\langle \frac{\partial}{\partial t}u\left(t\right),\phi_{p}\right\rangle \right)^{2},\label{eq:100}
\end{equation}

\begin{equation}
\rho_{2}\left(\beta,\lambda_{p},s,t,\epsilon\right)=T^{2}\sum_{p\ge1}\int_{0}^{t}\left(\Psi\left(\beta,\lambda_{p},s,t\right)-\frac{e^{\sqrt{\lambda_{p}}\left(s-t\right)}}{2\sqrt{\lambda_{p}}}\right)^{2}\left|\left\langle f\left(s,u\left(s\right)\right)-f\left(s,u^{\epsilon}\left(s\right)\right),\phi_{p}\right\rangle \right|^{2}ds.\label{eq:101}
\end{equation}

Next, we shall estimate these terms (\ref{eq:100})-(\ref{eq:101})
as follows.

\begin{equation}
\rho_{1}\le\frac{\beta^{2}}{4}\left(\frac{\beta}{T}\right)^{\frac{-2t}{T}}\left(\ln\left(\frac{T}{\beta}\right)\right)^{\frac{-2t}{T}}\sum_{p\ge1}e^{2\sqrt{\lambda_{p}}\left(T-t\right)}\left(\sqrt{\lambda_{p}}\left\langle u\left(t\right),\phi_{p}\right\rangle +\left\langle \frac{\partial}{\partial t}u\left(t\right),\phi_{p}\right\rangle \right)^{2},\label{eq:102}
\end{equation}

\begin{eqnarray}
\rho_{2} & \le & T^{2}\sum_{p\ge1}\int_{0}^{t}\left(2\Psi^{2}\left(\beta,\lambda_{p},s,t\right)+\frac{1}{2\lambda_{p}}\right)\left|\left\langle f\left(s,u\left(s\right)\right)-f\left(s,u^{\epsilon}\left(s\right)\right),\phi_{p}\right\rangle \right|^{2}ds\nonumber \\
 & \le & \frac{T^{2}}{\lambda_{1}}\sum_{p\ge1}\int_{0}^{t}\left(\frac{\beta}{T}\right)^{\frac{2s-2t}{T}}\left(\ln\left(\frac{T}{\beta}\right)\right)^{\frac{2s-2t}{T}}\left|\left\langle f\left(s,u\left(s\right)\right)-f\left(s,u^{\epsilon}\left(s\right)\right),\phi_{p}\right\rangle \right|^{2}ds\nonumber \\
 & \le & \frac{T^{2}K^{2}}{\lambda_{1}}\int_{0}^{t}\left(\frac{\beta}{T}\right)^{\frac{2s-2t}{T}}\left(\ln\left(\frac{T}{\beta}\right)\right)^{\frac{2s-2t}{T}}\left\Vert u\left(s\right)-u^{\epsilon}\left(s\right)\right\Vert ^{2}ds.\label{eq:103}
\end{eqnarray}

Combining (\ref{eq:99}) and (\ref{eq:102})-(\ref{eq:103}), we have

\begin{eqnarray}
\left\Vert u\left(t\right)-u^{\epsilon}\left(t\right)\right\Vert ^{2} & \le & \frac{\beta^{2}}{4}\left(\frac{\beta}{T}\right)^{\frac{-2t}{T}}\left(\ln\left(\frac{T}{\beta}\right)\right)^{\frac{-2t}{T}}P\nonumber \\
 &  & +\frac{T^{2}K^{2}}{\lambda_{1}}\int_{0}^{t}\left(\frac{\beta}{T}\right)^{\frac{2s-2t}{T}}\left(\ln\left(\frac{T}{\beta}\right)\right)^{\frac{2s-2t}{T}}\left\Vert u\left(s\right)-u^{\epsilon}\left(s\right)\right\Vert ^{2}ds,\label{eq:104}
\end{eqnarray}

where $P$ is defined as in (\ref{eq:23}).

Multiplying both sides of (\ref{eq:104}) by ${\displaystyle \left(\frac{\beta}{T}\right)^{\frac{2t}{T}}\left(\ln\left(\frac{T}{\beta}\right)\right)^{\frac{2t}{T}}}$,
it yields

\begin{equation}
{\displaystyle \left(\frac{\beta}{T}\right)^{\frac{2t}{T}}\left(\ln\left(\frac{T}{\beta}\right)\right)^{\frac{2t}{T}}}\left\Vert u\left(t\right)-u^{\epsilon}\left(t\right)\right\Vert ^{2}\le\beta^{2}P+\frac{T^{2}K^{2}}{\lambda_{1}}\int_{0}^{t}\left(\frac{\beta}{T}\right)^{\frac{2s}{T}}\left(\ln\left(\frac{T}{\beta}\right)\right)^{\frac{2s}{T}}\left\Vert u\left(s\right)-u^{\epsilon}\left(s\right)\right\Vert ^{2}ds.\label{eq:105}
\end{equation}

Applying Gronwall's inequality to (\ref{eq:105}), we conclude that

\begin{equation}
{\displaystyle \left(\frac{\beta}{T}\right)^{\frac{2t}{T}}\left(\ln\left(\frac{T}{\beta}\right)\right)^{\frac{2t}{T}}}\left\Vert u\left(t\right)-u^{\epsilon}\left(t\right)\right\Vert ^{2}\le e^{\frac{T^{2}K^{2}t}{2\lambda_{1}}}P\beta^{2},
\end{equation}

implies the estimate (\ref{eq:34}). $\square$

\subsection*{Proof of Theorem \ref{thm:2} and Theorem \ref{thm:8}.}

Substituting $\beta=\epsilon^{m}$ into the estimates in four lemmas
\ref{lem:3}-\ref{lem:6} and using triangle inequality, it is straightforward
to conclude the whole desired results Theorem \ref{thm:2}. Similarly,
substituting $\beta=\epsilon^{m}$ into the estimates in three lemmas
\ref{lem:9} and \ref{lem:11}-\ref{lem:12} and using triangle inequality
yield the estimate (\ref{eq:25}). Moreover, the uniqueness result
in Lemma \ref{lem:10} implies the uniqueness of $v^{\epsilon}$ mentioned
in Theorem \ref{thm:8}. $\square$
\end{document}